\documentclass[11pt]{article}

\usepackage{latexsym}
\usepackage{amsfonts}
\usepackage{amsmath}

\usepackage[dvips]{color}

\newtheorem{thrm}{Theorem}[section]
\newtheorem{lemma}[thrm]{Lemma}
\newtheorem{prop}[thrm]{Proposition}

\newtheorem{remark}[thrm]{Remark}

\numberwithin{equation}{section}

\usepackage{enumerate}

\usepackage{amssymb}
\usepackage{mathtools}
\mathtoolsset{showonlyrefs}
\usepackage{mathrsfs}
\usepackage{comment}
\usepackage{hyperref}
\usepackage{xcolor}

\def\E{\mathbb{E} }
\def\P{\mathbb{P} }

\def\R{\mathbb{R} }
\def\N{\mathbb{N} }

\makeatletter
\begin{document}
\allowdisplaybreaks

\title{\Large \bf{
Asymptotic expansion for  additive measure of branching Brownian motion
}\footnote{The research of this project is supported
     by the National Key R\&D Program of China (No. 2020YFA0712900).}}
\author{ \bf  Haojie Hou \hspace{1mm}\hspace{1mm}
Yan-Xia Ren\footnote{The research of this author is supported by NSFC (Grant Nos. 12071011 and 12231002) and
The Fundamental Research Funds for the Central Universities, Peking University LMEQF.
 } \hspace{1mm}\hspace{1mm} and \hspace{1mm}\hspace{1mm}
Renming Song\thanks{Research supported in part by a grant from the Simons
Foundation
(\#960480, Renming Song).}
\hspace{1mm} }
\date{}
\maketitle

\begin{abstract}
	
Let $N(t)$ be the collection of particles
	alive at time $t$ in a branching Brownian motion in $\R^d$, and for $u\in N(t)$, let $\mathbf{X}_u(t)$ be the position of particle $u$ at time $t$.
	For $\theta\in \R^d$, we define the additive measures of the branching Brownian motion by
	$$
		\mu_t^\theta (\mathrm{d}\mathbf{x}):= e^{-(1+\frac{\Vert\theta\Vert^2}{2})t}\sum_{u\in N(t)} e^{-\theta \cdot \mathbf{X}_u(t)} \delta_{\left(\mathbf{X}_u(t)+\theta t\right)}(\mathrm{d}\mathbf{x}).
	$$
	In this paper, under some conditions on the offspring distribution, we give
	asymptotic expansions of arbitrary order for
	$\mu_t^\theta ((\mathbf{a}, \mathbf{b}])$ and $\mu_t^\theta ((-\infty,  \mathbf{a}])$ for $\theta\in \R^d$ with $\Vert \theta \Vert <\sqrt{2}$.
	These expansions sharpen the asymptotic results of
	Asmussen and Kaplan (1976) and Kang (1999),
	and are analogs of the expansions in
	Gao and Liu (2021) and R\'{e}v\'{e}sz,  Rosen and Shi (2005)
	for branching Wiener processes (a particular class of branching random walks) corresponding to
	$\theta=\mathbf{0}$.

\end{abstract}

\medskip


\medskip


 \section{Introduction and main results}

\subsection{Introduction}

A branching random walk in $\R^d$ is a discrete-time Markov process which can be defined as follows:
at time $0$, there is a particle at $\mathbf{0}\in \R^d$; at time $1$, this particle is replaced by a random number of particle distributed according to a point process $\mathcal{L}$;
at time $2$, each individual,  of generation 1, if located at $\mathbf{x}\in \R^d$, is replaced by
a point process $\mathbf{x}+ \mathcal{L}_\mathbf{x}$, where $\mathcal{L}_\mathbf{x}$ is
an independent copy of $\mathcal{L}$.
This procedure goes on.  We use $Z_n$ to denote the point process formed by the positions of the particles of generation $n$.

Biggins \cite{Biggins92} studied the $L^p$ convergence of the additive martingale
\[
W_n(\theta):= \frac{1}{m(\theta)^n} \int e^{-\theta \cdot \mathbf{x}}Z_n(\mathrm{d} \mathbf{x}),
\]
where $m(\theta):= \mathbb{E}\left(\int e^{-\theta \cdot \mathbf{x}}Z_1(\mathrm{d} \mathbf{x})\right)$.
He used the $L^p$ convergence of the additive martingale to study the asymptotic behavior of
$Z_n(n{\bf c}+I)$ for fixed $\bf{c}$ and bounded interval $I$.
To describe Biggins' result, we introduce the following additive measure
$\mu_n^{Z, \theta}$ of the branching random walk, which is a shifted version of the measure introduced before
Theorem 4 in \cite{Biggins92}:
\begin{equation}\label{Additive-measure-Z}
	\mu_n^{Z, \theta}(A):= m(\theta)^{-n}  \int e^{-\theta \cdot \mathbf{y}} 1_{A}(\mathbf{y}-\mathbf{c}_\theta n)Z_n (\mathrm{d}\mathbf{y}),\quad A\in {\cal B}(\mathbb{R}^d),
\end{equation}
with $(\mathbf{c}_\theta)_i:=m(\theta)^{-1}  \mathbb{E}\left( \int  x_i e^{-\theta\cdot\bf{x}} Z_1 (\mathrm{d} \mathbf{x}) \right)$.
\cite[Theroem 4]{Biggins92} implies  that,
in the weak disorder regime (i.e., $-\log m(\theta)< - \theta \cdot \triangledown m(\theta) / m(\theta)$),
if there exists
$\gamma>1$
such that $\mathbb{E} (W_1(\theta)^\gamma)<\infty$, then
for $\mathbf{x}\in \R^d$ and $h>0$, as $n\to\infty$,
\begin{align}\label{Additive-measure-Z'}
	n^{d/2}\mu_n^{Z, \theta} (\mathbf{x}+ I_h) \longrightarrow \frac{(2h)^dW_\infty (\theta)}{(2\pi \mbox{det}(\Sigma_\theta))^{d/2} },
	\quad \mbox{a.\/s.}
\end{align}
where $W_\infty(\theta):= \lim_{n\to\infty} W_n (\theta), I_h= [-h,h]^d$
and
\[
(\Sigma_\theta)_{i,j}=
m(\theta)^{-1}\mathbb{E}\Big(\int (x_i- (\mathbf{c}_\theta)_i)(x_j- (\mathbf{c}_\theta)_j) e^{-\theta\cdot\bf{x}}Z_1(\mathrm{d} \mathbf{x})\Big),\quad i,j\in \{1,...,d\}.
\]
In the case $d=1$,  Pain proved that, see  \cite[(1.14)]{Pain18}, in the weak disorder regime, if there exists $\gamma>1$ such that $\mathbb{E} (W_1(\theta)^\gamma)<\infty$, then
for any $b\in \R$, as $n\to\infty$,
\begin{align}\label{Supercritical-CLT'}
	\mu_n^{Z, \theta} \left((-\infty, b\Sigma_\theta \sqrt{n}]\right) 	\to W_\infty (\theta) \Phi(b) \quad \text{ in  probability,}
\end{align}
where $\Phi(b):= \frac{1}{\sqrt{2\pi}}\int_{-\infty}^b e^{-z^2/2}\mathrm{d}z$.

For the case $\theta=\mathbf{0}$, there are many further asymptotic results.
In the case when $d=1$ and the point process  $\mathcal{L}$ is given by  $\mathcal{L}= \sum_{i=1}^{B} \delta_{X_i}$, where $X_i$
are iid with common   distribution $G$ and $B$ is an independent $\N$-valued random variable with  $\mathbb{P}(B=k)=p_k$ and $\mu:=\sum_kkp_k>1$,
Asmussen and Kaplan \cite{AK76-1, AK76-2} proved that if  $G$ has mean 0, variance 1
and $\sum_{k=2}^\infty k(\log k)^{1+\varepsilon} p_k <\infty$ for some $\varepsilon>0$, then
conditioned on survival, for any $b\in \mathbb{R}$,
\begin{align}\label{CLT}
	\mu^{Z, 0}_n
	((-\infty, b\sqrt{n}]) \stackrel{n\to\infty}{\longrightarrow}
	W_\infty(0) \Phi(b),\quad \mbox{a.s.}
\end{align}
They also proved that if $G$ has finite 3rd moment and $\sum_{k=2}^\infty k(\log k)^{3/2+\varepsilon} p_k <\infty$ for some $\varepsilon>0$,  then, for any $a<b\in \mathbb{R}$, conditioned on survival,
\begin{align}\label{LLT}
	\sqrt{2\pi n} \mu^{Z, 0}_n
	([a,b])\stackrel{n\to\infty}{\longrightarrow}
	(b-a)W_\infty(0),\quad \mbox{a.s.}
\end{align}
Gao and Liu \cite{GL16-2} gave first and second order expansions of $\mu^{Z, 0}_n((-\infty, b\sqrt{n}])$. A third order expansion was proved by Gao and Liu \cite{GL16-1, GL18}, where branching random walks in  (time) random environment were studied.
They also conjectured  the form of
asymptotic expansion of arbitrary order for $\mu^{Z, 0}_n((-\infty, b\sqrt{n}])$.
For general branching random walks, results similar to  \eqref{CLT} and \eqref{LLT} were proved in Biggins \cite{Biggins90}.

When the point process $\mathcal{L}$ is given by $\mathcal{L}= \sum_{i=1}^B \delta_{\mathbf{X}_i}$ where $\mathbf{X}_1,\mathbf{X}_2,...$ are
independent $d$-dimensional standard normal random variables and
$B$ is an independent $\N$-valued random variable with  $\mathbb{P}(B=k)=p_k$ and $\mu:=\sum_kkp_k>1$,
$Z_n$ is called a supercritical branching Wiener process.
R\'{e}v\'{e}sz \cite{RP94} first proved the analogs of \eqref{CLT} and \eqref{LLT} for branching Wiener processes,  then Chen \cite{Chen01} studied the corresponding convergence rates.
Gao and Liu \cite{GL21} proved that, for each $m\in \mathbb{N}$, when $\sum_{k=1}^\infty k(\log k)^{1+\lambda}p_k<\infty$ for some  $\lambda > 3\max \{ (m+1), dm  \}$, there exist random variables $\{V_{\mathbf{a}}, |\mathbf{a}| \leq m\}$ such that for each $\mathbf{t}\in \mathbb{R}^d$,
\[
\frac{1}{\mu^n} Z_n((-\infty, \mathbf{t}\sqrt{n}]) = \Phi_d(\mathbf{t}) V_{\mathbf{0}} +\sum_{\ell=1}^m \frac{(-1)^\ell}{n^{\ell/2}} \sum_{|\mathbf{a}| =\ell} \frac{D^\mathbf{a} \Phi_d(\mathbf{t})}{\mathbf{a}!}V_{\mathbf{a}}+ o(n^{-m/2}),\quad
\mbox{a.s.}
\]
where for $\mathbf{a}=(a_1,...,a_d)$, $|\mathbf{a}|= a_1+...+a_d$, $\mathbf{a}!= a_1 ! \cdots a_d!$, $\Phi_d(\mathbf{t})$ is the distribution function of
a $d$-dimensional standard normal random vector
and
$D^\mathbf{a} \Phi_d(\mathbf{t}):= \partial_{t_1}^{a_1}...\partial_{t_d}^{a_d}\Phi_d(\mathbf{t}).$
For the local limit theorem \eqref{LLT}, R\'{e}v\'{e}sz, Rosen and Shi \cite{RRS05} proved that, when $\sum_{k=1}^\infty k^2p_k<\infty$, for any bounded Borel set $A\subset \mathbb{R}^d$,
\begin{align}\label{Asymptotic-LIL}
	(2\pi n)^{d/2} \frac{1}{\mu^n}Z_n(A)= \sum_{\ell=0}^m\frac{(-1)^\ell}{(2n)^\ell}\sum_{|\mathbf{a}| =\ell} \frac{1}{\mathbf{a}!} \sum_{\mathbf{b}\leq 2\mathbf{a}} C_{2\mathbf{a}}^{\mathbf{b}}(-1)^{|\mathbf{b}|} M_{\mathbf{b}}(A) V_{2\mathbf{a}-\mathbf{b}} +o(n^{-m}),\quad
	\mbox{a.s.},
\end{align}
where $\mathbf{b}\leq 2\mathbf{a}$ means that $b_i\leq 2a_i$ for all $1\leq i\leq d$,
$C_{2\mathbf{a}}^{\mathbf{b}} := C_{2a_1}^{b_1}\cdots C_{2a_d}^{b_d}$ and $M_{\mathbf{b}}(A): = \int_A x_1^{b_1}\cdots x_d^{b_d}\mathrm{d}x_1...\mathrm{d}x_d$.

For the lattice case,
analogs of \eqref{CLT} and \eqref{LLT} can be found in \cite{Gao17, GK17}, and an  asymptotic expansion  similar
to \eqref{Asymptotic-LIL} for $Z_n(\{k\})$ was given by Gr\"{u}bel and Kabluchko\cite{GK17}.

In this paper, we are concerned with branching Brownian motions in $\R^d$.
A branching Brownian motion in $\R^d$ is a continuous-time Markov process defined as follows: initially there is a particle at $\mathbf{0}\in \mathbb{R}^d$, it moves according to a $d$-dimensional standard Brownian motion and its lifetime is an exponential random variable of parameter 1, independent of the spatial motion.
At the end of its lifetime, it produces $k$ offspring with probability $p_k$ for  $k\in \mathbb{N}$ and the offspring
move independently according to a $d$-dimensional standard Brownian motion from the death location of their parent,
and repeat their parent's behavior independently.
This procedure goes on.
We will use $\P$ to denote the law of branching Brownian motion and $\E$ to denote the corresponding expectation.
Without loss of generality, we assume that
\[
\sum_{k=0}^\infty kp_k =2.
\]
Let $N(t)$ be the set of particles alive at time $t$ and for $u\in N(t)$, we use $\mathbf{X}_u(t)$ to denote the position of particle $u$ at time $t$.  Define
\[
Z_t:= \sum_{u\in N(t)} \delta_{\mathbf{X}_u(t)}.
\]
For $\theta=(\theta_1,...,\theta_d) \in \mathbb{R}^d$,
\begin{align}
	W_t(\theta):= e^{-(1+\frac{\Vert\theta\Vert^2}{2})t}\sum_{u\in N(t)} e^{-\theta \cdot \mathbf{X}_u(t)}
\end{align}
is a non-negative martingale and is called the additive martingale of the branching Brownian motion.
When $\theta$ is the zero vector, $W_t(\theta)$ reduces to $e^{-t}Z_t(\R^d)$.
It is well-known that (for $d=1$, see Kyprianou\cite{Ky04}),
for each $\theta \in \R^d$, $W_t(\theta)$
converges to a non-trivial limit $W_\infty(\theta)$ if and only if $\Vert \theta\Vert<\sqrt{2}$ and
\begin{align}\label{Weaker-LlogL-moment}
	\sum_{k=1}^\infty k (\log k)p_k <\infty.
\end{align}

From now on, we will only consider $\theta \in \R^d$ with $\|\theta\|<\sqrt{2}$.
For any set $A\subset \R$ and $a\in \R$, we use $|A|$ to denote the Lebesgue measure of $A$ and $a A:= \{ax: x\in A\}.$
Asmussen and Kaplan \cite[Part 5]{AK76-2} proved that when $d=1$, under the assumption $\sum_{k=1}^\infty k^2 p_k<\infty$, for any Borel set $B$ with $|\partial B|=0$, as $t\to\infty$,
\begin{align}\label{Continuous-CLT}
	e^{-t}Z_t\left(\sqrt{t}B\right)\longrightarrow \ \frac{W_\infty(0)}{\sqrt{2\pi}}\int_B e^{-z^2/2}\mathrm{d}z, \quad
	\P\mbox{-a.s.}
\end{align}
and that for any bounded Borel set $B$ with $|\partial B|=0$, as $t\to\infty$,
\begin{align}\label{Continuous-LLT}
	\sqrt{2\pi t}e^{-t}Z_t\left( B\right)\longrightarrow \ |B| W_\infty(0), \quad
	\P\mbox{-a.s.}
\end{align}
Kang \cite[Theorem 1]{Kang99} weakened the moment condition and proved that \eqref{Continuous-CLT} holds with $B= (-\infty, b]$ under condition \eqref{Weaker-LlogL-moment}.

Similar to \eqref{Additive-measure-Z}, we define
the additive measure $\mu_t^\theta$ of branching Brownian motion as
\begin{align}\label{Additive-measure}
	\mu_t^\theta (\mathrm{d}\mathbf{x}):= e^{-(1+\frac{\Vert\theta\Vert^2}{2})t}\sum_{u\in N(t)} e^{-\theta \cdot \mathbf{X}_u(t)} \delta_{\left(\mathbf{X}_u(t)+\theta t\right)}(\mathrm{d}\mathbf{x}).
\end{align}
The aim of this paper is to prove asymptotic expansions of arbitrary order
for $\mu^\theta_t$ for $\theta\in \R^d$ with $\|\theta\|<\sqrt2$,
see Theorems \ref{thm1} and \ref{thm2} below.
These expansions sharpen the asymptotic results of  \cite[Part 5]{AK76-2, Kang99} mentioned above.
The asymptotic expansions of  \cite{GL21, RRS05} are for the additive measure $\mu^{Z, 0}_n$ of
branching Wiener processes, while the asymptotic expansions of  Theorems \ref{thm1} and \ref{thm2} are for
the additive measure $\mu^\theta_t$ of branching Brownian motions with $\theta$ not necessarily $\mathbf{0}$.

One might expect that the asymptotic expansions for branching Wiener processes, when considered along $\left\{t_n= n\delta, n\in \N \right\}$, can be used to used to get the expansions of this paper by letting $\delta\to 0$. However, it seems that this idea does not work due to two different reasons.  One of the reasons is that values along $\left\{n\delta, n\in \N \right\}$
are not good enough
to control the behavior between the time intervals $[t_n, t_{n+1}]$. Another reason is that
$\{Z_{n\delta}: n\in \N\}$ is not  a branching Wiener process since in
$Z_{\delta}=\sum_{u\in N(\delta)} \delta_{\mathbf{X}_u(\delta)}$,  for $u, v\in N(\delta), u\neq v$, $ \mathbf{X}_u(\delta)$ and $\mathbf{X}_v(\delta)$ are not independent.

\subsection{Notation}

We list here some notation that will be used repeatedly below.
Throughout this paper, $\N=\{0, 1, \dots\}$.
Recall that $N(t)$ is the set of the particles alive at time $t$ and that for $u\in N(t)$, $\mathbf{X}_u(t)$ is the position of $u$.
For $u\in N(t)$,
we use $d_u$ and $O_u$ to denote the death time and the offspring number of $u$ respectively.  For $v$ and $u$, we will use $v<u$ to denote that $v$ is an ancestor of $u$. The notation $v\leq u$ means that $v=u$ or $v<u$.

For $\mathbf{a}=(a_1,...,a_d)\in \R^d$,  define $(\mathbf{a})_j:= a_j$ and   $(-\infty, \mathbf{a}]:= (-\infty, a_1]\times \cdots \times (-\infty, a_d]$.
For $\mathbf{a}, \mathbf{b}\in \R^d$,
we use $\mathbf{a}< \mathbf{b}$ ($\mathbf{a}\leq \mathbf{b}$) to denote that $(\mathbf{a})_j< (\mathbf{b})_j $ ($(\mathbf{a})_j< (\mathbf{b})_j $ ) for all $1\leq j\leq d$.
For $\mathbf{a}, \mathbf{b}\in \R^d$
with $\mathbf{a}< \mathbf{b}$,
define $(\mathbf{a}, \mathbf{b}]:= (a_1, b_1]\times\cdots \times (a_d, b_d]$. The definition of $[\mathbf{a}, \mathbf{b}]$ is similar. For $\mathbf{k}=(k_1,...,k_d)\in \N^d$, set $|\mathbf{k}|:=k_1+...+k_d$ and $\mathbf{k}!:= k_1!\cdots k_d!$.
For a  function $f$ on $\R^d$, $\mathbf{x}\in \R^d$ and $\mathbf{k}\in \N^d$, let
$D^{\mathbf{k}}f(\mathbf{x}):= \partial_{x_1}^{k_1}...  \partial_{x_d}^{k_d}f(\mathbf{x})$.
We also use the notation $\phi(y):= \frac{1}{\sqrt{2\pi}} e^{-y^2/2}$ and $\Phi_d(\mathbf{x}):=\prod_{j=1}^d \int_{-\infty}^{x_j} \phi(z)\mathrm{d}z$. Sometimes we write $\Phi(y)$ for $\Phi_1(y)$.

\subsection{Main results}

We will assume that
\begin{align}\label{LlogL}
	\sum_{k=1}^\infty k(\log k)^{1+\lambda}p_k < 	\infty
\end{align}
for appropriate $\lambda>0$.
Let $H_k$ be the $k$-th order Hermite polynomial: $H_0(x):=1$ and for $k\geq 1$,
\[
H_k(x):= \sum_{j=0}^{[k/2]} \frac{k! (-1)^j}{2^j j! (k-2j)!} x^{k-2j}.
\]
It is well known that
if $\{ (B_t)_{t\geq 0}, \Pi_0\}$ is a standard Brownian motion, then for any $k\ge 0$,
$\{t^{k/2} H_k(B_t/\sqrt{t}), \sigma(B_s: s\leq t),\ \Pi_0\}$ is a martingale. Now for
$\mathbf{k}\in \N^d$ and $\theta\in \R^d$ with $\|\theta\|<\sqrt{2}$, we define
\begin{align}\label{Martingale}
	M_t^{(\mathbf{k},\theta)}:= e^{-(1+\frac{\Vert\theta\Vert^2}{2})t}\sum_{u\in N(t)} e^{-\theta \cdot \mathbf{X}_u(t)} t^{|\mathbf{k}|/2} \prod_{j=1}^d H_{k_j}\Big(\frac{(\mathbf{X}_u(t))_j+\theta_j t}{\sqrt{t}}\Big),
	\quad t\geq 0.
\end{align}
Note that $M_t^{({\bf 0}, \theta)}= W_t(\theta)$. We will prove in
Proposition \ref{Convergence-rate-Martingale} below that
if \eqref{LlogL} holds for $\lambda$ large enough, $M_t^{({\bf k},\theta)}$ will converges almost surely and  in $L^1$ to a limit $M_\infty^{({\bf k},\theta)}$.
Here are the main results of this paper:

\begin{thrm}\label{thm1}
	Suppose $\theta\in \R^d$ with $\Vert\theta\Vert<\sqrt{2}$.
	For any given $m\in \mathbb{N}$, if \eqref{LlogL} holds for some $\lambda> \max\left\{ 3m+8, d(3m+5)\right\}$,
	then for any $\mathbf{b}\in \mathbb{R}^d$, $\mathbb{P}$-almost surely, as $s\to\infty$,
	\begin{align}\label{Result-1}
		\mu_{s}^\theta \left((-\infty, \mathbf{b}\sqrt{s}]\right)  &= \sum_{\mathbf{k}: |\mathbf{k}|\leq m} \frac{(-1)^{|\mathbf{k}|}}{\mathbf{k}!} \frac{1}{s^{|\mathbf{k}|/2}} D^\mathbf{k} \Phi_d(\mathbf{b}) M_{\infty}^{(\mathbf{k},\theta)}+ o(s^{-m/2})\\
		& = \sum_{\ell=0}^m
		\frac{(-1)^\ell}{s^{\ell/2}}
		\sum_{\mathbf{k}: |\mathbf{k}|= \ell} \frac{D^\mathbf{k} \Phi_d(\mathbf{b})}{\mathbf{k}!}M_{\infty}^{(\mathbf{k},\theta)}+ o(s^{-m/2}).
	\end{align}
\end{thrm}

\begin{thrm}\label{thm2}
	Suppose $\theta\in \R^d$ with $\Vert\theta\Vert<\sqrt{2}$.
	For any given $m\in \mathbb{N}$, if \eqref{LlogL} holds for some $\lambda> \max\{d(3m+5), 3m+3d+8\}$,
	then for any
	$\mathbf{a}, \mathbf{b}\in \mathbb{R}^d$ with $\mathbf{a}< \mathbf{b}$,
	$\mathbb{P}$-almost surely, as $s\to\infty$,
	\begin{align}
		&s^{d/2}\mu_s^\theta
		((\mathbf{a}, \mathbf{b}])  \nonumber\\
		& = \sum_{\ell=0}^m \frac{1}{s^{\ell/2}} \sum_{j=0}^\ell (-1)^j\sum_{\mathbf{k}: |\mathbf{k}|=j}\frac{M_{\infty}^{(\mathbf{k},\theta)}}{\mathbf{k}!}
		\sum_{\mathbf{i}: |\mathbf{i}|=\ell-j}
		\frac{ D^{\mathbf{k}+\mathbf{i}+\mathbf{1}}\Phi_d(\mathbf{0})}{\mathbf{i}!}
		\int_{[\mathbf{a},\mathbf{b}]}
		\prod_{j=1}^d z_{j}^{i_j}\mathrm{d}z_1...\mathrm{d}z_d +o(s^{-m/2}),
	\end{align}
	where $\mathbf{1}:=(1,...,1)$.
\end{thrm}

\begin{remark}
	Note that we only dealt with the case that the branching rate is 1 and the mean number of offspring is 2 in the two theorems above. In the general case when the branching rate is $\beta>0$ and  the mean number of offspring is $\mu>1$, one can use the same argument to prove the following counterpart of Theorem \ref{thm1}: Suppose $\theta\in \R^d$ with $\Vert\theta\Vert<\sqrt{2\beta(\mu-1)}$.
	For any given $m\in \mathbb{N}$, if \eqref{LlogL} holds for some $\lambda> \max\left\{ 3m+8, d(3m+5)\right\}$,
	then for any $\mathbf{b}\in \mathbb{R}^d$, $\mathbb{P}$-almost surely, as $s\to\infty$,
	\begin{align}
		\mu_{s}^\theta \left((-\infty, \mathbf{b}\sqrt{s}]\right) &:= 	
		e^{-(\beta(\mu-1)+\frac{\Vert\theta\Vert^2}{2})t}
		\sum_{u\in N(t)} e^{-\theta \cdot \mathbf{X}_u(t)} 1_{ (-\infty, \mathbf{b}\sqrt{s}]}  \left( \mathbf{X}_u(t)+\theta t \right)
		\\ & = \sum_{\ell=0}^m
		\frac{(-1)^\ell}{s^{\ell/2}}
		\sum_{\mathbf{k}: |\mathbf{k}|= \ell} \frac{D^\mathbf{k} \Phi_d(\mathbf{b})}{\mathbf{k}!}M_{\infty}^{(\mathbf{k},\theta)}+ o(s^{-m/2}),
	\end{align}
	with $M_\infty^{(\mathbf{k},\theta)}$ given by
	\begin{align}\label{General-Martingale}
		M_\infty^{(\mathbf{k},\theta)}:=\lim_{t\to\infty}
		e^{-(\beta(\mu-1)+\frac{\Vert\theta\Vert^2}{2})t}
		\sum_{u\in N(t)} e^{-\theta \cdot \mathbf{X}_u(t)} t^{|\mathbf{k}|/2} \prod_{j=1}^d H_{k_j}\Big(\frac{(\mathbf{X}_u(t))_j+\theta_j t}{\sqrt{t}}\Big).
	\end{align}

	In the general case, the counterpart of Theorem \ref{thm2} is as follows:
	Suppose $\theta\in \R^d$ with $\Vert\theta\Vert<\sqrt{2\beta(\mu-1)}$.
	For any given $m\in \mathbb{N}$, if \eqref{LlogL} holds for some $\lambda> \max\{d(3m+5), 3m+3d+8\}$,
	then for any
	$\mathbf{a}, \mathbf{b}\in \mathbb{R}^d$ with $\mathbf{a}< \mathbf{b}$,
	$\mathbb{P}$-almost surely, as $s\to\infty$,
	\begin{align}
		&s^{d/2}\mu_s^\theta
		((\mathbf{a}, \mathbf{b}])  =
		e^{-(\beta(\mu-1)+\frac{\Vert\theta\Vert^2}{2})t}
		\sum_{u\in N(t)} e^{-\theta \cdot \mathbf{X}_u(t)} 1_{ (\mathbf{a}, \mathbf{b}]}  \left( \mathbf{X}_u(t)+\theta t \right)  \nonumber\\
		& = \sum_{\ell=0}^m \frac{1}{s^{\ell/2}} \sum_{j=0}^\ell (-1)^j\sum_{\mathbf{k}: |\mathbf{k}|=j}\frac{M_{\infty}^{(\mathbf{k},\theta)}}{\mathbf{k}!}
		\sum_{\mathbf{i}: |\mathbf{i}|=\ell-j}
		\frac{ D^{\mathbf{k}+\mathbf{i}+\mathbf{1}}\Phi_d(\mathbf{0})}{\mathbf{i}!}
		\int_{[\mathbf{a},\mathbf{b}]}
		\prod_{j=1}^d z_{j}^{i_j}\mathrm{d}z_1...\mathrm{d}z_d +o(s^{-m/2}),
	\end{align}
	with $M_\infty^{(\mathbf{k},\theta)}$ given in \eqref{General-Martingale}.
\end{remark}

\begin{remark}
	One could also consider asymptotic expansions for the additive measure $\mu^{Z, \theta}_n$ for branching random walks.
	Using the   tools established in \cite{GL16-1}, it is possible to get fixed order expansions. However, getting
	asymptotic expansions of arbitrary order
	may be difficult.
\end{remark}

We end this section with a few words about the strategy of the proofs and the organization of the paper.
In Section \ref{pre}, we introduce the spine decomposition and gather some useful facts.
We also study
the convergence rate of the martingales $M_t^{({\bf k},\theta)}$ and moments
of the additive martingale $W_t(\theta)$.
In Section \ref{Main}, we prove Theorems \ref{thm1} and \ref{thm2}.
To prove Theorem \ref{thm1}, we choose a sequence of discrete time
$r_n=n^{1/\kappa}$ for some $\kappa>1$.
To control the behavior of particles alive in $(r_n, r_{n+1})$, we need $r_{n+1}-r_n\to0$. This is the reason we do not choose $r_n=n\delta$.
We prove in Lemma \ref{lemma1} that
$  \mu_{r_n}^\theta \left((-\infty, \mathbf{b}\sqrt{r_n}]\right) \approx \mathbb{E} \left[ \mu_{r_n}^\theta \left((-\infty, \mathbf{b}\sqrt{r_n}]\right) \big| \mathcal{F}_{\sqrt{r_n}}\right]$,
where $\mathcal{F}_t$ is the $\sigma$-field generated by the branching Brownian motion
up to time $t$.
To deal with $s\in (r_n, r_{n+1})$, we adapt some ideas from from \cite[ Lemma 8]{AK76-2} and \cite[paragraph below (13)]{Kang99}.
We prove in Lemma \ref{lemma2} that, for $s\in (r_n, r_{n+1})$,
$\mu_{s}^\theta \left((-\infty, \mathbf{b}\sqrt{s}]\right) \approx \mathbb{E} \left[ \mu_{s}^\theta \left((-\infty, \mathbf{b}\sqrt{s}]\right) \big| \mathcal{F}_{\sqrt{r_n}}\right] $.
We complete the proof of  Theorem \ref{thm1}  by
using a series of identities proved in \cite{GL21}.
The proof of Theorem \ref{thm2} is similar.

\section{Preliminaries}\label{pre}

\subsection{Spine decomposition}

Define
\begin{align}\label{Change-of-measure}
	\frac{\mathrm{d}\mathbb{P}^{-\theta}}{\mathrm{d} \mathbb{P}} \bigg|_{\mathcal{F}_t}:= W_t(\theta).
\end{align}
Then under $\mathbb{P}^{-\theta}$,
the evolution of our branching Brownian motion can be described as follows (spine decomposition) (see \cite{Ky04} for the case $d=1$ or see \cite{RS20} for a more general case):

(i) there is an initial marked particle at $\mathbf{0} \in \mathbb{R}^d$ which moves
according to the law of $\{\mathbf{B}_t -\theta t, \Pi_\mathbf{0}\}$,
where $\{\mathbf{B}_t, \Pi_\mathbf{0}\}$ is a $d$-dimensional standard Brownian motion;

(ii) the branching rate of this marked particle is $2$;

(iii) when the marked particle dies at site $\mathbf{y}$,
it gives birth to $\widehat{L}$ children with $\mathbb{P}^{-\theta}(\widehat{L}=k)=kp_k/2$;

(iv)  one of these children is uniformly selected and marked, and
the marked child evolves as its parent independently
and the other children evolve independently with law $\mathbb{P}_\mathbf{y}$, where $\mathbb{P}_\mathbf{y}$ denotes the law of a branching Brownian motion starting at $\mathbf{y}$.

Let $d_i$ be the $i$-th splitting time of the spine and $O_i$ be the number of children produced by the spine at time $d_i$.  According to the spine decomposition, it is easy to see that $\{d_i: i\geq 1\}$ are the atoms for a Poisson point process with rate 2, $\{O_i: i\geq 1\}$ are iid with common law $\widehat{L}$ given by $\mathbb{P}^{-\theta}(\widehat{L}=k)= kp_k/2$, and that
$\left\{d_i: i\geq 1 \right\}$ and $\{O_i: i\geq 1\}$ and $\mathbf{X}_\xi$ are independent. This fact will be used repeatedly.

We use $\xi_t$  and $\mathbf{X}_\xi(t)$ to denote the marked particle at time $t$ and the position of this marked particle respectively.
By \cite[Theorem 2.11]{RS20}, we have that, for $u\in N(t)$,
\begin{align}\label{Many-to-one}
	\mathbb{P}^{-\theta}\left(\xi_t =u \big| \mathcal{F}_t\right) =
	\frac{e^{-\theta \cdot \mathbf{X}_u(t)}}{\sum_{u\in N(t)} e^{-\theta \cdot \mathbf{X}_u(t)}}=\frac{e^{-(1+\frac{\Vert\theta\Vert^2}{2})t } e^{-\theta \cdot \mathbf{X}_u(t)}}{W_t(\theta)}.
\end{align}
Using \eqref{Many-to-one},
we can  get the following many-to-one formula.

\begin{lemma}\label{General-many-to-one}
	For any $t>0$ and $u\in N(t)$, let $H(u, t)$ be a  non-negative $\mathcal{F}_t$-measurable  random variable. Then
	\[
	\E \Big(\sum_{u\in N(t)} H(u,t) \Big) =e^{(1+\frac{\Vert\theta\Vert^2}{2})t}\mathbb{E}^{-\theta}\left(e^{\theta \cdot \mathbf{X}_\xi(t)}H(\xi_t, t)\right).
	\]
\end{lemma}
\textbf{Proof: } Combining \eqref{Change-of-measure} and \eqref{Many-to-one}, we get
\begin{align}
	&	\E \left(\sum_{u\in N(t)} H(u,t) \right) = \E^{-\theta} \left(\sum_{u\in N(t)} \frac{H(u,t)}{W_t(\theta)} \right) \nonumber\\
	&= \E^{-\theta}\left(\sum_{u\in N(t)}  H(u,t)  	 e^{(1+\frac{\Vert\theta\Vert^2}{2})t } e^{\theta \cdot \mathbf{X}_u(t)}  \mathbb{P}^{-\theta}\left(\xi_t =u \big| \mathcal{F}_t\right)\right)\nonumber\\
	& = e^{(1+\frac{\Vert\theta\Vert^2}{2})t }  \E^{-\theta}\left(     \mathbb{E}^{-\theta}\left(  \sum_{u\in N(t)} 1_{ \left\{\xi_t =u \right\}} H(u,t)  	  e^{\theta \cdot \mathbf{X}_u(t)}  \big| \mathcal{F}_t\right)\right)\nonumber\\
	& =  e^{(1+\frac{\Vert\theta\Vert^2}{2})t } \E^{-\theta}\left(    H(\xi_t,t)  	 e^{\theta \cdot \mathbf{X}_\xi(t)}   \sum_{u\in N(t)} 1_{ \left\{\xi_t =u \right\}} \right) = e^{(1+\frac{\Vert\theta\Vert^2}{2})t}\mathbb{E}^{-\theta}\left(e^{\theta \cdot \mathbf{X}_\xi(t)}H(\xi_t, t)\right).
\end{align}
\hfill$\Box$

\subsection{Some useful facts}

In this subsection, we gather some useful  facts that will be used later.

\begin{lemma}\label{Useful-Ineq-1}
	(i) Let $\ell \in [1,2]$ be a fixed constant. Then for any
	finite family of independent centered random variables  $\left\{X_i : i=1, \dots, n\right\}$ with
	${\rm E}|X_i|^\ell <\infty$ for all $i=1, \dots, n$,
	it holds that
	\[
	{\rm E}
	\big|\sum^n_{i=1} X_i\big|^\ell	\le 2 \sum^n_{i=1}
	{\rm E}
	|X_i|^\ell.
	\]
	(ii) For any $\ell\in [1,2]$ and any random variable $X$ with
	${\rm E}|X|^2 <\infty$,
	\[
	{\rm E}	\left|X-{\rm E} X	\right|^\ell \lesssim 	{\rm E}	|X|^\ell \leq (	{\rm E}	X^2)^{\ell /2}.
	\]
\end{lemma}
\textbf{Proof: }
	For (i), see  \cite[Theorem 2]{BE65}. (ii) follows easily from Jensen's inequality.
\hfill$\Box$

\begin{lemma}\label{Asymptotic-expansion}
	For any $\rho\in (0,1), b,x \in \mathbb{R}$, it holds that
	\[
	\Phi\left(\frac{b-\rho x}{\sqrt{1-\rho^2}}\right) = \Phi (b)- \phi(b) \sum_{k=1}^\infty \frac{\rho^k}{k!} H_{k-1}(b) H_k(x).
	\]
\end{lemma}
\textbf{Proof: }
	See \cite[Lemma 4.2.]{GL21}.
\hfill$\Box$
\bigskip

To prove Theorem \ref{thm1},  we will define $r_n:= n^{\frac{1}{\kappa}}$ for some $\kappa>1$.
For $s\in [r_n, r_{n+1})$, applying Lemma \ref{Asymptotic-expansion} with $\rho= \sqrt{\sqrt{r_n} /s}$ and $x=r_n^{-1/4}y$, we  get that for any $b, y\in \R$,
\begin{align}\label{Asymptotic-expansion-2}
	\Phi\Big(\frac{b\sqrt{s}- y}{\sqrt{s-\sqrt{r_n}}}\Big) = \Phi (b)- \phi(b) \sum_{k=1}^\infty \frac{1}{k!} \frac{1}{s^{k/2}} H_{k-1}(b) \Big(r_n^{k/4} H_k\Big(\frac{y}{r_n^{1/4}}\Big)\Big).
\end{align}

Recall that (see \cite[(4.1)]{GL21}) for any $k\geq 1$ and $x\in \mathbb{R}$,
\begin{align}\label{Upper-H_k-2}
	\left|H_k(x)\right|\leq 2\sqrt{k!}e^{x^2/4}.
\end{align}

\begin{lemma}\label{lemma3}
	For a given $m\in \mathbb{N}$, let $\kappa= m+3$ and $r_n = n^{1/\kappa}$.
	Let $K>0$ be a fixed constant and $J$ be an integer such that $J > 2m+ K\kappa $. For any $b, y \in \mathbb{R}$
	and $s\in [r_n, r_{n+1})$, it holds that
	\[
	\Phi\Big(\frac{b\sqrt{s}- y}{\sqrt{s-\sqrt{r_n}}}\Big) = \Phi (b)- \phi(b) \sum_{k=1}^J \frac{1}{k!} \frac{1}{s^{k/2}} H_{k-1}(b) \Big(r_n^{k/4} H_k\Big(\frac{y}{r_n^{1/4}}\Big)\Big) +
	\varepsilon_{m,y,b, s},
	\]
	and that
	\[
	\sup\Big\{s^{m/2}	\left|\varepsilon_{m,y,b,s}\right|: \ s\in[r_n, r_{n+1}), \ |y| \leq \sqrt{K\sqrt{r_n}\log n}, b\in \R \Big\}
	\stackrel{n\to\infty}{\longrightarrow}0.
	\]
\end{lemma}
\textbf{Proof: }
	It follows from \eqref{Upper-H_k-2} that there exists a constant $C$
	such that for all $b\in \R$, $n\geq 2$, $s\in [r_n, r_{n+1})$, $|y|\leq \sqrt{K \sqrt{r_n} \log n}$ and $k\geq m$,
	\begin{align}\label{Uniform-Upper-Bound}
		s^{m/2}\frac{1}{k!}\frac{1}{s^{k/2}} \left(\phi(b)\left|H_{k-1}(b)\right| \right) \Big|r_n^{k/4} H_k\Big(\frac{y}{r_n^{1/4}}\Big)\Big|  \leq C
		\frac{n^{k/(4\kappa)}}{n^{(k-m)/(2\kappa)}} n^{K/4}.
	\end{align}
	Combining this with \eqref{Asymptotic-expansion-2}, we get that  for $J> 2m+ K\kappa $, $n\geq 2$, $s\in [r_n, r_{n+1})$,
	$b\in \R$ and $|y|\leq \sqrt{K \sqrt{r_n} \log n}$,
	\begin{align}
		s^{m/2}\left|\varepsilon_{m,y,b,s}\right| \leq C
		\sum_{k=J+1}^\infty n^{-(k-2m- K\kappa)/(4\kappa) } \lesssim
		n^{-(J+1-2m- K\kappa)/(4\kappa) }.
	\end{align}
	Thus the  assertions of the lemma are valid.
\hfill$\Box$
\bigskip

Now we give a result of similar flavor which will be used to prove Theorem \ref{thm2}.
Taking derivative with respect to $b$ in Lemma \ref{Asymptotic-expansion},
and using the fact that
\begin{align}\label{K-th-derivative}
	\frac{\mathrm{d}^k}{\mathrm{d} b^k} \Phi(b)=(-1)^{k-1} H_{k-1}(b)\phi(b),
\end{align}
we get that

\begin{align}\label{step_23}
	\frac{1}{\sqrt{1-\rho^2}} \phi \Big( \frac{b-\rho x}{\sqrt{1-\rho^2}}\Big) = \phi(b) + \phi(b) \sum_{k=1}^\infty \frac{\rho^k}{k!}H_k(b)H_k(x).
\end{align}
Now letting $\rho= \sqrt{\sqrt{r_n}/s}$, $b=z/\sqrt{s}$ and
$x= r_n^{-1/4}y$ in \eqref{step_23}, we get that for any $z, y\in \R$,
\begin{align}\label{step_25}
	\frac{\sqrt{s}}{\sqrt{s-\sqrt{r_n}}}\phi\Big( \frac{z-y}{\sqrt{s-\sqrt{r_n}}}\Big) = \phi\Big(\frac{z}{\sqrt{s}}\Big) \Big(1+ \sum_{k=1}^\infty \frac{1}{k!} \frac{1}{s^{k/2}} H_k\Big(\frac{z}{\sqrt{s}}\Big)  r_n^{k/4}H_k\Big(\frac{y}{r_n^{1/4}}\Big) \Big).
\end{align}

The proof of the following result is similar to that of Lemma \ref{lemma3}, we omit the details.

\begin{lemma}\label{lemma4}
	For a given $m\in \mathbb{N}$, let $\kappa= m+3$ and $r_n = n^{1/\kappa}$. Let $K>0$ be a fixed constant and $J$ be an integer such that $J > 2m+ K\kappa $. For any $a< b\in \mathbb{R}$, $y, z\in \R$ and $s\in [r_n, r_{n+1})$, it holds that
	\begin{align}
	& \frac{\sqrt{s}}{\sqrt{s-\sqrt{r_n}}}\phi\Big( \frac{z-y}{\sqrt{s-\sqrt{r_n}}}\Big)\\
	& = \phi\Big(\frac{z}{\sqrt{s}}\Big) \Big(1+ \sum_{k=1}^J \frac{1}{k!} \frac{1}{s^{k/2}} H_k\Big(\frac{z}{\sqrt{s}}\Big)  r_n^{k/4}H_k\Big(\frac{y}{r_n^{1/4}}\Big) \Big)+ \varepsilon_{m,y, z, s},
	\end{align}
	and that
	\[
	\sup\Big\{  s^{m/2} \left|\varepsilon_{m,y,z, s}\right|: \ s\in [r_n, r_{n+1}), \ z\in [a,b],\ |y| \leq \sqrt{K\sqrt{r_n}\log n} \Big\}
	\stackrel{n\to\infty}{\longrightarrow}0.
	\]
\end{lemma}

\subsection{Convergence rate for the martingales}

\begin{prop}
	\label{Convergence-rate-Martingale}
	For any $\theta\in \R^d$ with $\Vert\theta\Vert<\sqrt{2}$ and $\mathbf{k}\in \N^d$,
	$\{M_t^{(\mathbf{k}, \theta )}, t\geq 0; \,\mathbb{P}\}$ is a martingale.
	If \eqref{LlogL} holds for  some $\lambda > |\mathbf{k}|/2$, then $M_t^{(\mathbf{k},\theta)}$ converges to a limit $M_\infty^{(\mathbf{k},\theta)}$ $\mathbb{P}$-a.s. and in $L^1$.
	Moreover,  for any $\eta\in (0, \lambda - |\mathbf{k}|/2)$,
	as $t\to\infty$,
	\[
	M_t^{(\mathbf{k},\theta)} - M_\infty^{(\mathbf{k},\theta)}=  o(t^{-(\lambda -|\mathbf{k}|/2)+\eta}),\quad 	
	\mathbb{P}\mbox{-a.s.}
	\]
\end{prop}
\textbf{Proof: }
By Lemma \ref{General-many-to-one}, it is easy to see that
$\{M_t^{(\mathbf{k}, \theta )}, t\geq 0; \,\mathbb{P}\}$ is a martingale.

Now we fix $\mathbf{k}\in \N^d$ and assume  \eqref{LlogL} holds for  some $\lambda > |\mathbf{k}|/2$.
We first look at the case when $t\to \infty$  along integers. Let $t=n\in \mathbb{N}$.
Recall that $N(n+1)$ is the set of particles alive at time $n+1$.
For $u\in N(n+1)$, define
$B_{n,u}$ to be the event that, for all $v<u$ with $d_v\in (n, n+1)$, it holds that $O_v \leq e^{c_0n}$,
where $c_0>0$ is a small constant to be determined later.
Set
\begin{align}
&M_{n+1}^{(\mathbf{k},\theta), B}:= e^{-(1+\frac{\Vert\theta\Vert^2}{2})(n+1)}\\
& \quad\quad\times \sum_{u\in N(n+1)} e^{-\theta \cdot \mathbf{X}_u(n+1)} (n+1)^{|\mathbf{k}|/2} \prod_{j=1}^d H_{k_j}\Big(\frac{\left(\mathbf{X}_u(n+1)\right)_j+\theta_j (n+1)}{\sqrt{n+1}}\Big)1_{B_{n,u}}.
\end{align}
Since $|H_k(x)|\lesssim |x|^k +1$ for all $x\in \mathbb{R}$ and $(|x|+|y|)^k \lesssim |x|^k +|y|^k$ for all $x, y\in \R$, we have
\begin{align}\label{Upper-H_k}
	(n+1)^{k/2} \Big|H_k \Big(\frac{x+z}{\sqrt{n+1}}\Big)\Big|\lesssim (|x|+|z|)^k + (n+1)^{k/2}\lesssim |x|^k + |z|^k +n^{k/2},
\end{align}
which implies that for all $j\in \{1,...,d\}$,
\begin{align}\label{Upper-H_k-3}
	& (n+1)^{k_j/2} \Big|H_{k_j}\Big(\frac{(\mathbf{X}_u(n+1))_j+\theta_j (n+1)}{\sqrt{n+1}}\Big) \Big|\nonumber\\
	& \lesssim |(\mathbf{X}_u(n))_j+\theta_j n|^{k_j}+ n^{k_j/2}+|(\mathbf{X}_u(n+1))_j-(\mathbf{X}_u(n))_j+\theta_j|^{k_j}.
\end{align}
Therefore,
\begin{align}
	&\left|M_{n+1}^{(\mathbf{k},\theta)} - M_{n+1}^{(\mathbf{k},\theta), B} \right| \leq e^{-(1+\frac{\Vert\theta\Vert^2}{2})(n+1)}\\
	&\quad\times \sum_{u\in N(n+1)} e^{-\theta \cdot \mathbf{X}_u(n+1)}  \Big| (n+1)^{|\mathbf{k}|/2} \prod_{j=1}^d H_{k_j}\Big(\frac{\left(\mathbf{X}_u(n+1)\right)_j+\theta_j (n+1)}{\sqrt{n+1}}\Big) \Big|1_{(B_{n,u})^c}\nonumber\\
	& \lesssim e^{-(1+\frac{\Vert\theta\Vert^2}{2})(n+1)}\sum_{u\in N(n)} e^{-\theta \cdot \mathbf{X}_u(n)}
	\sum_{v\in N(n+1): u\leq v}  e^{-\theta \cdot \left(\mathbf{X}_v(n+1)- \mathbf{X}_u(n)\right)} \nonumber\\
	& \quad \times \prod_{j=1}^d \left( |(\mathbf{X}_u(n))_j+\theta_j n|^{k_j}+ n^{k_j/2}+|(\mathbf{X}_v(n+1))_j-(\mathbf{X}_u(n))_j+\theta_j|^{k_j}\right)1_{(B_{n,v})^c}.
\end{align}
By the branching property and the Markov property, we get that
\begin{align}\label{First-moment-1}
	&  \mathbb{E}\left(\left| M_{n+1}^{(\mathbf{k},\theta)} -  M_{n+1}^{(\mathbf{k},\theta), B}\right| \big| \mathcal{F}_n \right) \lesssim e^{-(1+\frac{\Vert\theta\Vert^2}{2})n}\sum_{u\in N(n)} e^{-\theta \cdot \mathbf{X}_u(n)} \nonumber
	\\ & \,\,\,  \times	\mathbb{E}
	\Big( e^{-(1+\frac{\Vert\theta\Vert^2}{2})}\sum_{v\in N(1)} e^{-\theta \cdot \mathbf{X}_v(1)} \prod_{j=1}^d\Big( \left|(\mathbf{X}_v(1))_j+ \theta_j\right|^{k_j} +y_j\Big)1_{(D_{n,v})^c}  \Big)\Big|_{y_j=|(\mathbf{X}_u(n))_j+\theta_j n|^{k_j}+ n^{k_j/2}}
	\nonumber \\ &=: e^{-(1+\frac{\Vert\theta\Vert^2}{2})n}\sum_{u\in N(n)} e^{-\theta \cdot \mathbf{X}_u(n)} F\left(\mathbf{y}\right)\Big|_{y_j=|(\mathbf{X}_u(n))_j+
		\theta_j n|^{k_j}+ n^{k_j/2}},
\end{align}
where, for $v\in N(1)$, $D_{n,v}$ denotes the event that, for all $w<v$, it holds that $O_w\leq e^{c_0 n}$.
Recall that $d_i$ is the $i$-th splitting time of the spine and $O_i$ is the number of children produced by the spine at time $d_i$.
Define $D_{n, \xi_1}$
to be the event that, for all $i$ with $d_i<1$, it holds that $O_i\leq e^{c_0 n}$.
By Lemma \ref{General-many-to-one},
\begin{align}\label{First-moment-2}
	F(\mathbf{y})
	&= \mathbb{E}^{-\theta}\Big(\prod_{j=1}^d\left(\left|(\mathbf{X}_\xi(1))_j+\theta_j\right|^{k_j}+y_j\Big)
	1_{(D_{n,\xi_1})^c}\right).
\end{align}
Using the independence of $\{d_i: i\geq 1\}$, $\{O_i: i\geq 1\}$ and $\mathbf{X}_\xi$, we have that $D_{n,\xi_1}$ is independent of $\mathbf{X}_\xi$,  which implies that for $y_j\ge 1$,
\begin{align}\label{First-moment-3}
	& F(\mathbf{y})
	= \prod_{j=1}^d \left( \Pi_0 (|B_1|^{k_j})+y_j\right) \mathbb{P}^{-\theta}
	\left( D_{n,\xi_1}^c \right)
	\leq
	\prod_{j=1}^d \left( \Pi_0 (|B_1|^{k_j})+y_j\right) \mathbb{E}^{-\theta}
	\big( \sum_{i: d_i \leq 1} 1_{\{O_i> e^{c_0 n}\}}  \big) \nonumber\\
	& \lesssim \big(\prod_{j=1}^d y_j\big) \mathbb{P}^{-\theta}(\widehat{L}> e^{c_0 n})
	\lesssim \big(\prod_{j=1}^d y_j\big) \mathbb{P}^{-\theta}(\widehat{L}> e^{c_0 n}) \lesssim \frac{\prod_{j=1}^d y_j}{n^{1+\lambda}},
\end{align}
here in the first equality, we also used the fact that
$\{\mathbf{X}_\xi(t)+\theta t, \mathbb{P}^{-\theta}\}$ is a $d$-dimensional standard Brownian
motion, and in the last inequality, we used $\mathbb{E}^{-\theta} (\log_+^{1+\lambda}\widehat{L})<\infty$ (which follows from
\eqref{LlogL}).
By \eqref{First-moment-1} and \eqref{First-moment-3},  we have
\begin{align}
	& \mathbb{E}\left(\left| M_{n+1}^{(\mathbf{k},\theta)} -  M_{n+1}^{(\mathbf{k},\theta), B}\right| \big| \mathcal{F}_n \right) \\
	&\lesssim e^{-(1+\frac{\Vert\theta\Vert^2}{2})n}\sum_{u\in N(n)} e^{-\theta \cdot \mathbf{X}_u(n)}  \frac{ \prod_{j=1}^d \left(|(\mathbf{X}_u(n))_j+\theta_j n|^{k_j} +n^{k_j/2}\right) }{ n^{1+\lambda}}.
\end{align}
Taking expectation with respect to $\mathbb{P}$, by Lemma \ref{General-many-to-one}, we obtain that
\begin{align}\label{First-moment-4}
	& \mathbb{E}\left(\left|M_{n+1}^{(\mathbf{k},\theta)}-M_{n+1}^{(\mathbf{k},\theta), B}\right|\right) \lesssim  \frac{1}{n^{1+\lambda}} \prod_{j=1}^d\Pi_0 \left(|B_n|^{k_j} +n^{k_j/2}\right)  \lesssim \frac{n^{|\mathbf{k}|/2}}{n^{1+\lambda}}.
\end{align}
On the other hand, by the branching property,
\[
M_{n+1}^{(\mathbf{k},\theta), B}= e^{-(1+\frac{\Vert\theta\Vert^2}{2})n}\sum_{u\in N(n)} e^{-\theta \cdot \mathbf{X}_u(n)}  J_{n,u},
\]
where
\[
J_{n,u}:=
\sum_{v\in N(n+1): u\leq v}
e^{-\theta \cdot \left(\mathbf{X}_v(n+1)- \mathbf{X}_u(n)\right)}  (n+1)^{|\mathbf{k}|/2} \prod_{j=1}^d H_{k_j}\Big(\frac{(\mathbf{X}_v(n+1))_j+\theta_j (n+1)}{\sqrt{n+1}}\Big)1_{B_{n,v}}
\]
are independent given $\mathcal{F}_n$.
For any fixed $1< \ell < \min\{2/\Vert\theta\Vert^2, 2\}$,
Applying Lemma \ref{Useful-Ineq-1} (i) to the finite family $\{J_{n, u}-\mathbb{E}(J_{n,u}\big| \mathcal{F}_n): u\in N(n)\}$
and Lemma \ref{Useful-Ineq-1} (ii) to $J_{n,u}$,
together with \eqref{Upper-H_k-3}, we get that
\begin{align}\label{Second-moment-1}
	& \mathbb{E}\Big( \Big| M_{n+1}^{(\mathbf{k},\theta), B} - \mathbb{E} \Big( M_{n+1}^{(\mathbf{k},\theta), B}\big| \mathcal{F}_n \Big) \Big|^\ell \big| \mathcal{F}_n \Big)\lesssim e^{-\ell (1+\frac{\Vert \theta\Vert^2}{2})n}\sum_{u\in N(n)} e^{-\ell \theta \cdot \mathbf{X}_u(n)} \left(M_{n,u}\right)^{\ell/2},
\end{align}
where $M_{n,u}$ is given by
\begin{align}
	M_{n,u}:=  \mathbb{E}\Big( e^{-2(1+\frac{\Vert\theta\Vert^2}{2})}  \Big(\sum_{v\in N(1)} e^{-\theta \cdot \mathbf{X}_v(1)} S_v(\mathbf{y},1)1_{D_{n,v}}   \Big)^2 \Big)\Big|_{y_j= |(\mathbf{X}_u(n))_j+\theta_j n|^{k_j} +n^{k_j/2}}
\end{align}
with $S_v(\mathbf{y},r):= \prod_{j=1}^d\left( \left|(\mathbf{X}_v(r))_j+ \theta_j r\right|^{k_j} +y_j \right)$.
Set
\[
T_{n,u}:= S_u(\mathbf{y},1) 1_{D_{n,u}}e^{-(1+\frac{\Vert\theta\Vert^2}{2})}  \sum_{v\in N(1)} e^{-\theta \cdot \mathbf{X}_v(1)} S_v(\mathbf{y},1)1_{D_{n,v}}.
\]
By Lemma \ref{General-many-to-one},
we have
\begin{align}\label{Second-moment-2}
	&M_{n,u}= e^{-(1+\frac{\Vert\theta\Vert^2}{2})}   \mathbb{E}\Big(\sum_{u\in N(1)} e^{-\theta \cdot \mathbf{X}_u(1)} T_{n,u} \Big)\Big|_{y_j= |(\mathbf{X}_u(n))_j+\theta_j n|^{k_j} +n^{k_j/2}} \nonumber\\
	& =    \mathbb{E}^{-\theta}\left(T_{n,\xi_1} \right)\Big|_{y_j= |(\mathbf{X}_u(n))_j+\theta_j n|^{k_j} +n^{k_j/2}} \nonumber\\
	&= \mathbb{E}^{-\theta} \Big(
	1_{D_{n,\xi_1}}
	S_{\xi_1}(\mathbf{y},1) e^{-(1+\frac{\Vert\theta\Vert^2}{2})} \sum_{v\in N(1)} e^{-\theta \cdot \mathbf{X}_v(1)} S_v(\mathbf{y},1)1_{D_{n,v}} \Big)
	\Big|_{y_j= |(\mathbf{X}_u(n))_j+\theta_j n|^{k_j} +n^{k_j/2}}\nonumber\\
	& \leq \mathbb{E}^{-\theta} \Big(
	1_{D_{n,\xi_1}}
	S_{\xi_1}(\mathbf{y}, 1) e^{-(1+\frac{\Vert\theta\Vert^2}{2})} \sum_{v\in N(1)} e^{-\theta \cdot \mathbf{X}_v(1)} S_v(\mathbf{y},1) \Big)
	\Big|_{y_j= |(\mathbf{X}_u(n))_j+\theta_j n|^{k_j} +n^{k_j/2}}.
\end{align}
Conditioned on $\mathcal{G}:= \sigma(\mathbf{X}_\xi, d_i, O_i: i\geq 1)$,
by the branching property,  on the set $D_{n,\xi_1}$,
\begin{align}\label{Second-moment-3}
	& e^{-(1+\frac{\Vert \theta\Vert^2}{2})} 	\mathbb{E}^{-\theta} \Big(  \sum_{v\in N(1)} e^{-\theta \cdot \mathbf{X}_v(1)} S_v(\mathbf{y},1)\big| \mathcal{G} \Big) \nonumber \\
	& = \sum_{i: d_i \leq 1}  (O_i -1)e^{-(1+\frac{\Vert\theta\Vert^2}{2})d_i} e^{-\theta \cdot \mathbf{X }_\xi (d_i)} \Pi_{\mathbf{0}} \Big(\prod_{j=1}^d\left( \left|(\mathbf{B}_{1-d_i})_j+ z_j\right|^{k_j} +y_j \right)\Big)\Big|_{z_j=(\mathbf{X}_\xi(d_i))_j+\theta_j d_i}\nonumber\\
	& \lesssim e^{c_0 n} \sum_{i:d_i \leq 1} e^{-\theta \cdot \mathbf{X}_\xi (d_i)} \prod_{j=1}^d \left(y_j+ \left|(\mathbf{X}_\xi(d_i))_j+\theta_j d_i\right|^{k_j}\right) = e^{c_0 n} \sum_{i:d_i \leq 1} e^{-\theta \cdot \mathbf{X}_\xi (d_i)} S_{\xi_{d_i}}(\mathbf{y}, d_i) ,
\end{align}
where in the first equality we also used Lemma \ref{General-many-to-one}.
Plugging \eqref{Second-moment-3} into  \eqref{Second-moment-2},
noting that $\mathbf{X}_\xi$ and $d_i$ are independent, we get that
\begin{align}\label{Second-moment-4}
	& \mathbb{E}^{-\theta} \Big( 1_{D_{n,\xi_1}}S_{\xi_1}(\mathbf{y},1) e^{-(1+\frac{\Vert\theta\Vert^2}{2})} \sum_{v\in N(1)} e^{-\theta \cdot \mathbf{X}_v(1)} S_{v}(\mathbf{y},1)\big| \mathbf{X}_\xi   \Big)\nonumber\\
	& \lesssim 	S_{\xi_1}(\mathbf{y},1) 	e^{c_0 n} \mathbb{E}^{-\theta} \Big(\sum_{i:d_i \leq 1} e^{-\theta \cdot \mathbf{X}_\xi (d_i)}  S_{\xi_{d_i}}(\mathbf{y}, d_i)  \big| \mathbf{X}_\xi \Big)\nonumber\\
	& = 	 2S_{\xi_1}(\mathbf{y},1)
	e^{c_0 n} \int_0^1 e^{-\theta \cdot \mathbf{X}_\xi(s)} S_{\xi_{s}}(\mathbf{y}, s)   \mathrm{d}s \nonumber\\
	& \lesssim  e^{c_0 n}  \prod_{j=1}^d \left\{\left(y_j+\sup_{s<1} \left|(\mathbf{X}_\xi(s))_j+ \theta_j s\right|^{k_j}\right)^2 e^{\Vert\theta\Vert   \sup_{s<1} |(\mathbf{X}_\xi(s))_j+\theta_j s| }\right\}.
\end{align}
Since $\{\mathbf{X}_\xi (s)+\theta s, \mathbb{P}^{-\theta}\}$ is a $d$-dimensional standard Brownian motion,
combining with \eqref{Second-moment-2} and \eqref{Second-moment-4}, we conclude that
\begin{align}\label{Second-moment-5}
	& M_{n,u}\lesssim e^{c_0 n} \prod_{j=1}^d \Pi_\mathbf{0} \Big(\big(y_j+ \sup_{s<1}|(\mathbf{B}_s)_j|^{k_j}\Big)^2 e^{\Vert \theta\Vert \sup_{s<1} |(\mathbf{B}_s)_j|} \Big)\Big|_{y_j= |(\mathbf{X}_u(n))_j+\theta_j n|^{k_j} +n^{k_j/2}}\nonumber \\
	&\lesssim \prod_{j=1}^d\left(|(\mathbf{X}_u(n))_j+\theta_j n|^{k_j} +n^{k_j/2}\right)^2 e^{c_0 n}.
\end{align}
Plugging \eqref{Second-moment-5} into \eqref{Second-moment-1}, we conclude that
\begin{align}\label{Second-moment-6}
	& \mathbb{E}\Big( \left| M_{n+1}^{(\mathbf{k},\theta), B} - \mathbb{E} \big( M_{n+1}^{(\mathbf{k},\theta), B}\big| \mathcal{F}_n \big) \right|^\ell \big| \mathcal{F}_n \Big)\nonumber\\
	&\lesssim e^{c_0 \ell n/2}e^{-\ell (1+\frac{\Vert\theta\Vert^2}{2})n}\sum_{u\in N(n)} e^{-\ell \theta \cdot \mathbf{X}_u(n)} \prod_{j=1}^d\Big(|(\mathbf{X}_u(n))_j+\theta_j n|^{k_j} +n^{k_j/2}\Big)^\ell \nonumber\\
	& = e^{-\left((\ell -1)(1-\Vert\theta\Vert^2 \ell/2)-c_0 \ell /2\right)n }
	e^{-(1+\ell^2 \Vert\theta\Vert^2/2)n }\nonumber\\
	&\quad \times \sum_{u\in N(n)} e^{-\ell \theta \cdot \mathbf{X}_u(n)} \prod_{j=1}^d\left(|(\mathbf{X}_u(n))_j+\theta_j n|^{k_j} +n^{k_j/2}\right)^\ell.
\end{align}
Choose $c_0>0$ small so that $c_0 \ell /2< (\ell -1)(1-\Vert\theta\Vert^2\ell/2)$ and set $c_1:= (\ell -1)(1-\Vert\theta\Vert^2\ell/2) -c_0 \ell /2>0$. Taking expectation with respect to $\mathbb{P}$ in \eqref{Second-moment-6},
by Lemma \ref{General-many-to-one} with $\theta$ replaced to $\ell \theta$, we get that
\begin{align}\label{Second-moment-7}
	& \mathbb{E}\left( \left| M_{n+1}^{(\mathbf{k},\theta), B} - \mathbb{E} \left( M_{n+1}^{(\mathbf{k},\theta), B}\big| \mathcal{F}_n \right) \right|^{\ell} \right)\nonumber\\
	&\lesssim  e^{-c_1 n } \prod_{j=1}^d
	\Pi_\mathbf{0} \left(|(\mathbf{B}_n)_j - (\ell-1)\theta_j n |^{k_j} + n^{k_j/2}
	\right)^{\ell} \lesssim \left(\frac{n^{|\mathbf{k}|/2}}{n^{1+\lambda}}\right)^{\ell}.
\end{align}
Now combining \eqref{First-moment-4} and \eqref{Second-moment-7}, using the inequality:
\begin{align}\label{New-step_1}
	& {\rm E}\left(|X- {\rm E} (X| \mathcal{F}) |\right) \leq {\rm E} \left(|X-Y|\right)+ {\rm E}\left( |Y- {\rm E} (Y| \mathcal{F}) |\right) + {\rm E} \left(\left|{\rm E}\left(X-Y \big| \mathcal{F}\right) \right| \right)\nonumber\\
	& \leq 2 {\rm E} \left(|X-Y|\right) + {\rm E}\left( |Y- {\rm E} (Y| \mathcal{F}) |^\ell \right) ^{1/\ell},
\end{align}
and the fact that
$M_n^{(\mathbf{k},\theta)}= \mathbb{E} \left(M_{n+1}^{(\mathbf{k},\theta)}\big| \mathcal{F}_n\right)$,
we get that
\begin{align}\label{Moment-rate}
	& \mathbb{E}\left(\left|M_{n+1}^{(\mathbf{k},\theta)}-M_{n}^{(\mathbf{k},\theta)}\right|\right) \leq 2\mathbb{E}\left(\left|M_{n+1}^{(\mathbf{k},\theta)}-M_{n+1}^{(\mathbf{k},\theta), B}\right|\right)  \nonumber \\
	&\quad\quad+  \mathbb{E}\left( \left| M_{n+1}^{(\mathbf{k},\theta), B} - \mathbb{E} \left( M_{n+1}^{(\mathbf{k},\theta), B}\big| \mathcal{F}_n \right) \right|^{\ell} \right)^{1/\ell} \lesssim \frac{n^{|\mathbf{k}|/2}}{n^{1+\lambda}}.
\end{align}
Since $\lambda> |\mathbf{k}|/2$,  we have $\sum_{n=1}^\infty \mathbb{E}\left(\left|M_{n+1}^{(\mathbf{k},\theta)}-M_{n}^{(\mathbf{k},\theta)}\right|\right) <\infty$, which implies that
$M_n^{(\mathbf{k},\theta)}$ converges to a limit $M_\infty^{(\mathbf{k},\theta)}$ $\mathbb{P}$-almost surely and in $L^1$. Therefore, $M_n^{(\mathbf{k},\theta)}= \mathbb{E}\left(M_\infty^{(\mathbf{k},\theta)} \big| \mathcal{F}_n\right)$, $n\geq 1$.

For $s\in (n , n+1)$, $M_s^{(\mathbf{k},\theta)} =\mathbb{E}\left(M_{n+1}^{(\mathbf{k},\theta)}\big|\mathcal{F}_s\right)= \mathbb{E}\left(M_{\infty}^{(\mathbf{k},\theta)}\big| \mathcal{F}_s\right)$,
thus the second assertion of the proposition is valid.

Now we prove the last assertion of the proposition.
For any $\eta\in(0,\lambda -|\mathbf{k}|/2)$,  by \eqref{Moment-rate},
\begin{align}
	\sum_{n=1}^\infty n^{-|\mathbf{k}|/2 +\lambda -\eta } \mathbb{E}\left(\left|M_{n+1}^{(\mathbf{k},\theta)}-M_{n}^{(\mathbf{k},\theta)}\right|\right)  \lesssim \sum_{n=1}^\infty \frac{1}{n^{1+\eta}}<\infty,
\end{align}
which implies that
\begin{align}
	\sum_{n=1}^\infty n^{\lambda -|\mathbf{k}|/2 -\eta} \left(M_{n+1}^{(\mathbf{k},\theta)} - M_n^{(\mathbf{k},\theta)} \right)\ \mbox{converges }
	\mbox{a.s.}
\end{align}
Thus
$n^{\lambda -|\mathbf{k}|/2 -\eta}\left(M_n^{(\mathbf{k},\theta)}-M_\infty^{(\mathbf{k},\theta)}\right)\stackrel{n\to\infty}{\longrightarrow}0,\ \mathbb{P}$-a.s. (see for example \cite[Lemma 2]{Asmussen76}).
For $s\in[n, n+1]$, by Doob's inequality, for any $\varepsilon>0$,
\begin{align}
	&\sum_{n=1}^\infty \mathbb{P}\Big( n^{-|\mathbf{k}|/2 +\lambda -\eta} \sup_{n\leq s\leq n+1} \left| M_s^{(\mathbf{k},\theta)} -M_n^{(\mathbf{k},\theta)}\right| > \varepsilon \Big)\nonumber\\
	&\leq \frac{1}{\varepsilon} \sum_{n=1}^\infty n^{-|\mathbf{k}|/2 +\lambda -\eta } \mathbb{E}\Big(\Big|M_{n+1}^{(\mathbf{k},\theta)}-M_{n}^{(\mathbf{k},\theta)}\Big|\Big) <\infty.
\end{align}
Therefore, $n^{-|\mathbf{k}|/2+\lambda-\eta} \sup_{n\leq s\leq n+1} \left| M_s^{(\mathbf{k},\theta)} -M_n^{(\mathbf{k},\theta)}\right|  \stackrel{n\to\infty}{\longrightarrow} 0, \mathbb{P}$-a.s.
Hence $\mathbb{P}$-almost surely,
\begin{align}
	&\sup_{n\leq s\leq n+1} s^{-|\mathbf{k}|/2 + \lambda -\eta } \left|M_s^{(\mathbf{k},\theta)} - M_\infty^{(\mathbf{k},\theta)} \right|\\
	&\leq (n+1)^{-|\mathbf{k}|/2 +\lambda -\eta} \sup_{n\leq s\leq n+1} \left| M_s^{(\mathbf{k},\theta)} -M_n^{(\mathbf{k},\theta)}\right| + (n+1)^{-|\mathbf{k}|/2 +\lambda -\eta} \left|M_n^{(\mathbf{k},\theta)}-M_\infty^{(\mathbf{k}, \theta)}\right|\\
&\quad\stackrel{n\to\infty}{\longrightarrow}0,
\end{align}
which completes the proof of the last assertion of the proposition.
\hfill$\Box$

\subsection{Moment estimate for the additive martingale}

In this subsection, we give an upper bound for $W_t(\theta)$ which will be used later.

\begin{lemma}\label{LlogL-moment}
	Suppose $\theta\in \R^d$ with $\Vert \theta\Vert<\sqrt{2}$. If \eqref{LlogL} holds for some $\lambda>0$, then there exists a constant $C_{\theta,\lambda}$ such that for all $t>0$,
	\[
	\mathbb{E}\left(\left(W_t(\theta)+1\right)\log^{1+\lambda}(W_t(\theta)+1)\right)\leq C_{\theta , \lambda} (t+1).
	\]
\end{lemma}
\textbf{Proof:}
Since $\mathbb{E} W_t(\theta) =1$, it suffices to prove
that there exists a constant $C_{\theta,\lambda}$ such that for all $t>0$,
$\mathbb{E}\left(W_t(\theta)\log_+^{1+\lambda}\left(W_t(\theta)\right)\right)\le C_{\theta , \lambda} (t+1)$.
By using an projection argument, we can easily reduce to the one dimensional case. So we will only deal with the case $d=1$.
By \eqref{Change-of-measure}, we have $	\mathbb{E}\left(\left(W_t(\theta)\right)\log_+^{1+\lambda}(W_t(\theta))\right) = \mathbb{E}^{-\theta}\left( \log_+^{1+\lambda}(W_t(\theta))\right)$. Using the spine decomposition, we have
\begin{align}
	W_t(\theta)= e^{-(1+\frac{\theta^2}{2})t} e^{-\theta X_\xi(t)} + \sum_{i: d_i\leq t} e^{-(1+\frac{\theta^2}{2})d_i} e^{-\theta X_\xi(d_i)} \sum_{j=1}^{O_i-1} W_{t-d_i}^{i,j},
\end{align}
here $d_i, O_i$ are the $i$-th fission time and the number of offspring of the spine at time $d_i$ respectively.
Given all the information $\mathcal{G}$ about the spine,
$(W_{t-d_i}^{i,j})_{j\geq 1}$ are independent
with the same law as  $W_{t-d_i}(\theta)$ under $\mathbb{P}$.

Using elementary analysis one can easily show
that there exists $A=A_\lambda>1$ such that for any $x,y>A$,
\begin{align}\label{Ineq-1}
	\log_+^{1+\lambda}(x+y)\leq \log_+^{1+\lambda} (x)+ \log_+^{1+\lambda} (y).
\end{align}
We set
\begin{align}
	K_1&:= e^{-(1+\frac{\theta^2}{2})t} e^{-\theta X_\xi(t)} ,\\
	K_2&:= \sum_{i: d_i\leq t} e^{-(1+\frac{\theta^2}{2})d_i} e^{-\theta X_\xi(d_i)} \sum_{j=1}^{O_i-1} W_{t-d_i}^{i,j} 1_{\big\{e^{-(1+\frac{\theta^2}{2})d_i} e^{-\theta X_\xi(d_i)} \sum_{j=1}^{O_i-1} W_{t-d_i}^{i,j} \leq A\big\}} \leq A  \sum_{i: d_i\leq t} 1,\\
	K_3&:=  \sum_{i: d_i\leq t} e^{-(1+\frac{\theta^2}{2})d_i} e^{-\theta X_\xi(d_i)} \sum_{j=1}^{O_i-1} W_{t-d_i}^{i,j} 1_{\big\{e^{-(1+\frac{\theta^2}{2})d_i} e^{-\theta X_\xi(d_i)} \sum_{j=1}^{O_i-1} W_{t-d_i}^{i,j} > A\big\}}.
\end{align}
Note that $\log_+^{1+\lambda}(x+y+z)\leq \log_+^{1+\lambda}(3x)+  \log_+^{1+\lambda}(3y)+ \log_+^{1+\lambda}(3z)$,  $\log_+^{1+\lambda}(xy)\leq (\log_+ x + \log_+ y)^{1+\lambda} \lesssim \log_+^{1+\lambda}(x) + \log_+^{1+\lambda}(y)$ and  $\log_+^{1+\lambda}(x) \lesssim x$. By \eqref{Ineq-1}, we have
\begin{align}\label{Ineq-4}
	&\log_+^{1+\lambda} \left(W_t(\theta)\right)= \log_+^{1+\lambda} \left(K_1+K_2+K_3\right) \leq \log_+^{1+\lambda}(3K_1)+  \log_+^{1+\lambda}(3K_2)+ \log_+^{1+\lambda}(3K_3)\nonumber \\
	& \lesssim 1+\log_+^{1+\lambda}(K_1) +\Big( \sum_{i: d_i\leq t} 1 \Big)+  \sum_{i: d_i\leq t} \log_+^{1+\lambda} \Big(e^{-(1+\frac{\theta^2}{2})d_i} e^{-\theta X_\xi(d_i)} \sum_{j=1}^{O_i-1} W_{t-d_i}^{i,j}\Big)\nonumber\\
	& \lesssim 1+\log_+^{1+\lambda}(K_1) +\Big( \sum_{i: d_i\leq t} 1+   \log_+^{1+\lambda}  \big(e^{-(1+\frac{\theta^2}{2})d_i} e^{-\theta X_\xi(d_i)} \big) \Big)\nonumber\\
	& \quad \quad +  \sum_{i: d_i\leq t} \log_+^{1+\lambda}  \Big(\sum_{j=1}^{O_i-1} W_{t-d_i}^{i,j}\Big).
\end{align}
Put $\gamma := 1-\theta^2/2>0$.
Recalling that
$\{ X_\xi(t)+\theta t, \mathbb{P}^{-\theta}\}$
is a standard Brownian motion and $\{d_i: i\geq 1\}$ are the atoms of a Poisson point process with rate $2$ independent of $X_\xi$, we have
\begin{align}\label{Ineq-5}
	& \mathbb{E}^{-\theta}\Big(\log_+^{1+\lambda}(K_1) +\big( \sum_{i: d_i\leq t} 1+   \log_+^{1+\lambda}  \big(e^{-(1+\frac{\theta^2}{2})d_i} e^{-\theta X_\xi(d_i)} \big) \big)\Big) \nonumber\\
	& = \mathbb{E}^{-\theta}\Big(\log_+^{1+\lambda}(K_1) +2\int_0^t  \big(1+   \log_+^{1+\lambda}  \big(e^{-(1+\frac{\theta^2}{2})s} e^{-\theta X_\xi(s)} \big) \mathrm{d} s\big) \Big)\nonumber\\
	&\lesssim \Pi_0 \Big( \left(-\theta B_t -\gamma t\right)_+^{1+\lambda} + \int_0^t \big(1+ \left(-\theta B_s -\gamma s\right)_+^{1+\lambda} \big) \mathrm{d} s\Big)\lesssim t+1,
\end{align}
where the last inequality follows from the following estimate:
\begin{align}
	&\Pi_0 \big(\left(-\theta B_s -\gamma s\right)_+^{1+\lambda}\big) = s^{(1+\lambda)/2}\Pi_0 \big(\left(|\theta| B_1 -\gamma \sqrt{s}\right)^{1+\lambda}1_{\{|\theta|B_1>\gamma \sqrt{s}\}}\big)  \nonumber\\
	& \leq \gamma^{1+\lambda} (s+1)^{1+\lambda} e^{-\gamma \sqrt{s}} \Pi_0\big(\left(|\theta| |B_1| +1 \right) e^{|\theta| B_1}\big)\lesssim 1.
\end{align}
For the last term on the right-hand side of \eqref{Ineq-4}, conditioned on $\{d_i , O_i: i\ge 1\}$, we get
\begin{align}\label{Ineq-6}
	&\mathbb{E}^{-\theta} \Big(\sum_{i: d_i\leq t} \log_+^{1+\lambda}  \big(\sum_{j=1}^{O_i-1} W_{t-d_i}^{i,j}\big)\big| d_i, O_i, i\geq 1\Big)\nonumber\\
	&\lesssim   \sum_{i: d_i\leq t} \log_+^{1+\lambda} (O_i-1) + \sum_{i: d_i\leq t} \mathbb{E}^{-\theta} \Big( \log_+^{1+\lambda}  \big(\max_{j\leq O_i-1} W_{t-d_i}^{i,j}\big)\big| d_i, O_i, i\geq 1\Big).
\end{align}
Note that
\begin{align}\label{Ineq-7}
	&\mathbb{E}^{-\theta} \Big( \log_+^{1+\lambda}  \big(\max_{j\leq O_i-1} W_{t-d_i}^{i,j}\big)\big| d_i, O_i, i\geq 1\Big) \nonumber\\
	&= (1+\lambda)\int_0^\infty y^\lambda  \mathbb{P}^{-\theta} \Big( \max_{j\leq O_i-1} W_{t-d_i}^{i,j}>e^y\bigg| d_i, O_i, i\geq 1\Big)\mathrm{d} y \nonumber\\
	& = (1+\lambda)\int_0^\infty y^\lambda \Big(1-\prod_{j=1}^{O_i-1} \big( 1- \mathbb{P}^{-\theta} \big(  W_{t-d_i}^{i,j}>e^y\big| d_i, O_i, i\geq 1\big) \big) \Big) \mathrm{d} y \nonumber\\
	& \lesssim \int_0^\infty y^\lambda \big(1-\left( 1- e^{-y} \right)^{O_i-1} \big) \mathrm{d} y,
\end{align}
where in the inequality we used Markov's inequality.
When $O_i-1 \geq e^{y/2}$ (which is equivalent to $y\leq 2\log (O_i-1)$),
we have $y^\lambda \left(1-\left( 1- e^{-y} \right)^{O_i-1} \right)  \lesssim \log^{\lambda} (O_i-1)$; when $O_i-1 < e^{y/2}$, by the inequality $(1-x)^n\geq 1-nx$, we get
\[
y^\lambda \big(1-\left( 1- e^{-y} \right)^{O_i-1} \big) \leq y^\lambda (O_i-1) e^{-y}\leq y^\lambda e^{-y/2}.
\]
Thus, by \eqref{Ineq-7},
\begin{align}
	&\mathbb{E}^{-\theta} \Big( \log_+^{1+\lambda}  \big(\max_{j\leq O_i-1} W_{t-d_i}^{i,j}\big)\big| d_i, O_i, i\geq 1\Big) \lesssim  \log_+^{1+\lambda}(O_i-1) +1.
\end{align}
Plugging this back to \eqref{Ineq-6} and taking expectation with respect to $\mathbb{P}^{-\theta}$, we conclude that
\begin{align}\label{Ineq-8}
	&\mathbb{E}^{-\theta} \Big(\sum_{i: d_i\leq t} \log_+^{1+\lambda}  \big(\sum_{j=1}^{O_i-1} W_{t-d_i}^{i,j}\big)\Big)\nonumber\\
&\lesssim \mathbb{E}^{-\theta}\Big(\sum_{i: d_i\leq t} \big(\log_+^{1+\lambda} (O_i-1) +1\big)\Big) \lesssim t+1.
\end{align}
Combining  \eqref{Ineq-4}, \eqref{Ineq-5} and \eqref{Ineq-8}, we get the desired result.
\hfill$\Box$

\section{Proof of the main results}\label{Main}

Let $\kappa >1$ be fixed. Define
\begin{align}
	r_n:= n^{\frac{1}{\kappa}}, \quad n\in \mathbb{N}
\end{align}

\begin{lemma}\label{lemma1}
	For any given $\alpha,\beta >0$ and
	$\delta\in (0,1]$,
	assume that \eqref{LlogL} holds for  $\lambda$ with $\lambda \delta -\alpha >\kappa(1+\beta)$.
	
	(i) For each $n$, let $ a_n\leq n^\beta$ and $\{Y_{n,u}: u\in N(r_n^\delta)\}$ be a family of random variables
	such that $\mathbb{E}\left(Y_{n,u}\big| \mathcal{F}_{r_n^\delta}\right) =0$, and conditioned on $\mathcal{F}_{r_{n}^\delta}$, $Y_{n,u}, u\in N(r_n^\delta)$, are independent.
	If $|Y_{n,u}| \leq W_{a_n}(\theta; u) +1$ for all $n$ and $u\in N(r_n^\delta)$,
	with $\left(W_{a_n}(\theta; u), \P(\cdot\big| \mathcal{F}_{r_n^\delta})\right)$ being a copy of $W_{a_n}(\theta)$, then
	\begin{align}
		r_n^\alpha e^{-(1+\frac{\Vert\theta\Vert^2}{2})r_n^\delta} \sum_{u\in N(r_n^\delta) } e^{-\theta \cdot \mathbf{X}_u(r_n^\delta)}Y_{n,u} \stackrel{n\to\infty}{\longrightarrow} 0,\quad \mbox{a.s.}
	\end{align}
	
	(ii) Consequently, if $\lambda \delta -\alpha > \kappa +1$, then for any
	sequence $\{A_n\}$ of Borel sets in $\R^d$,
	\begin{align}
		r_n^{\alpha}\left| \mu_{r_n}^\theta \left(A_n\right) - \mathbb{E} \left[ \mu_{r_n}^\theta \left(A_n\right) \big| \mathcal{F}_{r_n^\delta}\right]  \right| \stackrel{n\to\infty}{\longrightarrow} 0,\quad
		\mbox{a.s.}
	\end{align}
\end{lemma}
\textbf{Proof:}
(i) Define
\[
\overline{Y}_{n,u}:= Y_{n,u}1_{\{|Y_{n,u}|\leq e^{c_*r_n^\delta} \}},\quad Y_{n,u}' = \overline{Y}_{n,u}- \mathbb{E}\left( \overline{Y}_{n,u} \big| \mathcal{F}_{r_n^\delta} \right),
\]
where $c_*>0$ is a constant to be chosen later.
Then for any $\varepsilon>0$,
\begin{align}\label{step_1}
	& \mathbb{P}\Big( \Big| r_n^\alpha e^{-(1+\frac{\Vert \theta\Vert^2}{2})r_n^\delta} \sum_{u\in N(r_n^\delta) } e^{-\theta \cdot \mathbf{X}_u(r_n^\delta)}Y_{n,u}  \Big|  > \varepsilon \Big| \mathcal{F}_{r_n^\delta} \Big) \nonumber\\
	& \leq  \mathbb{P}\Big( \Big| r_n^\alpha e^{-(1+\frac{\Vert \theta\Vert^2}{2})r_n^\delta} \sum_{u\in N(r_n^\delta) } e^{-\theta \cdot \mathbf{X}_u(r_n^\delta)}\left(Y_{n,u} - \overline{Y}_{n,u} \right) \Big|  > \frac{\varepsilon}{3} \Big| \mathcal{F}_{r_n^\delta} \Big)\nonumber\\
	& \quad
	+ \mathbb{P}\Big(r_n^\alpha e^{-(1+\frac{\Vert \theta\Vert^2}{2})r_n^\delta} \Big| \sum_{u\in N(r_n^\delta) } e^{-\theta \cdot \mathbf{X}_u(r_n^\delta)}Y_{n,u} ' \Big| > \frac{\varepsilon}{3}\Big| \mathcal{F}_{r_n^\delta}\Big)\nonumber\\
	& \quad + 1_{\big\{ r_n^\alpha e^{-(1+\frac{\Vert\theta\Vert^2}{2})r_n^\delta} \left| \sum_{u\in N(r_n^\delta) } e^{-\theta \cdot \mathbf{X}_u(r_n^\delta)}\mathbb{E}\left( \overline{Y}_{n,u} \big| \mathcal{F}_{r_n^\delta} \right)\right| >
		\frac{\varepsilon}{3} \big\}}
	=:I+II+III.
\end{align}
Using the inequality
\[
|Y_{n,u}- \overline{Y}_{n,u}|= |Y_{n,u}|1_{\{|Y_{n,u}|> e^{c_* r_n^\delta}\}} \leq
\left(W_{a_n}(\theta;u)+1\right)1_{\left\{W_{a_n}(\theta;u)+1> e^{c_* r_n^\delta}\right\}}
\]
and Markov's inequality, we have
\begin{align}\label{step_3}
	I &
	\leq \frac{3}{\varepsilon} r_n^\alpha e^{-(1+\frac{\Vert \theta\Vert^2}{2})r_n^\delta} \sum_{u\in N(r_n^\delta) } e^{-\theta \cdot \mathbf{X}_u(r_n^\delta)}\mathbb{E}\left( \left|Y_{n,u} - \overline{Y}_{n,u} \right|    \Big| \mathcal{F}_{r_n^\delta}\right)\nonumber\\
	& \lesssim r_n^\alpha e^{-(1+\frac{\Vert \theta\Vert^2}{2})r_n^\delta} \sum_{u\in N(r_n^\delta) } e^{-\theta \cdot \mathbf{X}_u(r_n^\delta)}\mathbb{E}\left( \left(W_{a_n}(\theta) +1\right) 1_{\{ W_{a_n}(\theta) +1>e^{c_*r_n^\delta}\}}  \right)\nonumber\\
	& \leq   \frac{r_n^\alpha}{(c_*r_n^\delta)^{\lambda+1}} e^{-(1+\frac{\Vert \theta\Vert^2}{2})r_n^\delta} \sum_{u\in N(r_n^\delta) } e^{-\theta \cdot \mathbf{X}_u(r_n^\delta)}  \mathbb{E}\left( \left(W_{a_n}(\theta) +1\right)\log_+^{1+\lambda}  \left(W_{a_n}(\theta) +1\right) \right) \nonumber\\
	& \lesssim \frac{r_n^\alpha}{(c_*r_n^\delta)^{\lambda+1}} e^{-(1+\frac{\Vert \theta\Vert^2}{2})r_n^\delta} \sum_{u\in N(r_n^\delta) } e^{-\theta \cdot \mathbf{X}_u(r_n^\delta)} n^\beta,
\end{align}
where in the last inequality we used Lemma \ref{LlogL-moment}.
For $1< \ell < \min\{ 2/\Vert\theta\Vert^2, 2\}$,
set $b:= e^{(\ell-1)(1-\Vert\theta\Vert^2 \ell /2)/2}\in (1, e)$ and $c_*:= \ln b$.
Using Markov's inequality, Lemma \ref{Useful-Ineq-1},
the conditional independence of $Y_{n,u}'$, and the fact that $|\overline{Y}_{n,u}|^\ell \leq e^{c_* (\ell -1)r_n^\delta} |Y_{n,u}|\leq e^{c_* r_n^\delta} \left( W_{a_n}(\theta;u)+1\right),$
we have
\begin{align}\label{step_4}
	II &
	\leq  \frac{3^\ell}{\varepsilon^\ell } r_n^{\alpha \ell}e^{-\ell(1+\frac{\Vert \theta\Vert^2}{2})r_n^\delta}   \mathbb{E}\Big(\Big| \sum_{u\in N(r_n^\delta) } e^{-\theta \cdot \mathbf{X}_u(r_n^\delta)}Y_{n,u} ' \Big|^\ell \Big| \mathcal{F}_{r_n^\delta}\Big)\nonumber\\
	&\lesssim
	\frac{3^\ell}{\varepsilon^\ell }
	r_n^{\ell \alpha} e^{-\ell (1+\frac{\Vert\theta\Vert^2}{2})r_n^\delta} \sum_{u\in N(r_n^\delta) } e^{-\ell \theta\cdot  \mathbf{X}_u(r_n^\delta)}\mathbb{E}\left( \left|\overline{Y}_{n,u}  \right|^\ell \big| \mathcal{F}_{r_n^\delta}\right)\nonumber \\
	& \leq  \frac{3^\ell}{\varepsilon^\ell } 	r_n^{\ell \alpha} e^{-\ell (1+\frac{\Vert\theta\Vert^2}{2})r_n^\delta} \sum_{u\in N(r_n^\delta) } e^{-\ell \theta\cdot  \mathbf{X}_u(r_n^\delta)}  e^{c_* r_n^\delta}\mathbb{E}\left(  W_{a_n} (\theta;u)+1 \big| \mathcal{F}_{r_n^\delta}\right)    \nonumber \\
	& \lesssim	 r_n^{\ell \alpha} b^{-2r_n^\delta}e^{-(1+\frac{\ell^2 \Vert\theta\Vert^2}{2})r_n^\delta} b^{r_n^\delta} \sum_{u\in N(r_n^\delta) } e^{-\ell \theta\cdot \mathbf{X}_u(r_n^\delta)},
\end{align}
where in the last inequality,
we used the identities
$\mathbb{E}\left(  W_{a_n} (\theta;u)+1 \big| \mathcal{F}_{r_n^\delta}\right)= \E\left(W_{a_n}(\theta)+1\right)=2$, $e^{-\ell (1+\frac{\Vert\theta\Vert^2}{2})r_n^\delta} = b^{-2r_n^\delta}e^{-(1+\frac{\ell^2 \Vert\theta\Vert^2}{2})r_n^\delta} $, and $e^{c_*}=b$.
Therefore, by \eqref{step_4}, we get
\begin{align}\label{step_6}
	II   \lesssim r_n^{\ell \alpha}  b^{-r_n^\delta} W_{r_n^\delta}(\ell \theta).
\end{align}
Now taking expectation with respect to $\mathbb{P}$ in \eqref{step_1}, and using \eqref{step_3} and \eqref{step_6}, we get that
\begin{align}\label{step_7}
	&  \mathbb{P}\Big( \Big| r_n^\alpha e^{-(1+\frac{\Vert\theta\Vert^2}{2})r_n^\delta} \sum_{u\in N(r_n^\delta) } e^{-\theta\cdot \mathbf{X}_u(r_n^\delta)}Y_{n,u}  \Big|  > \varepsilon  \Big)  \nonumber\\
	& \lesssim
	\frac{r_n^\alpha n^\beta}{(r_n^\delta)^{\lambda+1}}
	+ r_n^{\ell \alpha}  b^{-r_n^\delta} +
	\mathbb{P} \left(r_n^\alpha e^{-(1+\frac{\Vert\theta\Vert^2}{2})r_n^\delta} \Big| \sum_{u\in N(r_n^\delta) } e^{-\theta\cdot \mathbf{X}_u(r_n^\delta)}\mathbb{E}\left( \overline{Y}_{n,u} \big| \mathcal{F}_{r_n^\delta} \right)\Big| >
	\frac{\varepsilon}{3}\right).
\end{align}
By Markov's inequality and  the fact that
$\mathbb{E}\left( \overline{Y}_{n,u} \big| \mathcal{F}_{r_n^\delta} \right) = -\mathbb{E}\left( Y_{n,u}1_{\{|Y_{n,u}| >
	e^{c_*r_n^\delta}\} }
\big| \mathcal{F}_{r_n^\delta} \right)$,
the third term in right-hand side of \eqref{step_7} is bounded from above by
\begin{align}\label{step_8}
	&
	\frac{3}{\varepsilon} r_n^\alpha
	e^{-(1+\frac{\Vert\theta\Vert^2}{2})r_n^\delta} \mathbb{E} \Big(\sum_{u\in N(r_n^\delta) } e^{-\theta\cdot \mathbf{X}_u(r_n^\delta)}\mathbb{E}\big( |Y_{n,u}|
	1_{ \big\{|Y_{n,u}| >	e^{c_* r_n^\delta} \big\}}
	\big| \mathcal{F}_{r_n^\delta} \big) \Big)\nonumber \\
	& \leq \frac{3}{\varepsilon}r_n^\alpha  (c_* r_n^\delta)^{-\lambda-1} \mathbb{E}\left(\left(W_{a_n}(\theta)+1\right)\log_+^{\lambda+1} (1+W_{a_n}(\theta))\right)\lesssim r_n^\alpha (r_n^\delta)^{-\lambda-1} n^\beta,
\end{align}
where in the last inequality we used Lemma \ref{LlogL-moment}.
Plugging the upper bound above into  \eqref{step_7} and recalling $r_n= n^{\frac{1}{\kappa}}$, we get
\begin{align}
	& \sum_{n=1}^\infty \mathbb{P}\Big( \Big| r_n^\alpha e^{-(1+\frac{\Vert\theta\Vert^2}{2})r_n^\delta} \sum_{u\in N(r_n^\delta) } e^{-\theta\cdot  \mathbf{X}_u(r_n^\delta)}Y_{n,u}  \Big|  > \varepsilon  \Big)  \\ &\lesssim \sum_{n=1}^\infty \Big(
	\frac{r_n^\alpha n^\beta }{(r_n^\delta)^{\lambda+1}}
	+ r_n^{\ell \alpha}  b^{-r_n^\delta}
	+ r_n^{-((\lambda+1) \delta -\alpha-\kappa\beta )}\Big),
\end{align}
which is summable since $\lambda \delta - \alpha > \kappa(1+\beta)$. This completes the proof of (i).

(ii)
By the Markov property
and Lemma \ref{General-many-to-one},
\begin{align}
	& \mathbb{E} \left[ \mu_{r_n}^\theta \left(A_n\right) \big| \mathcal{F}_{r_n^\delta}\right] \\
	&= e^{-(1+\frac{\Vert\theta\Vert^2}{2})r_n^\delta} \sum_{u\in N(r_n^\delta)} e^{-\theta \cdot \mathbf{X}_u(r_n^\delta)} \mathbb{P}^{-\theta}\left(\mathbf{X}_\xi(r_n -r_n^\delta)+ \theta r_n +\mathbf{y} \in A_n\right)\big|_{\mathbf{y}= \mathbf{X}_u(r_n^\delta)}.
\end{align}
Since $\{\mathbf{X}_\xi(t)+\theta t, \mathbb{P}^{-\theta}\}$ is a $d$-dimensional standard Brownian motion, we have
\begin{align}
	& \mathbb{E} \left[ \mu_{r_n}^\theta \left(A_n\right) \big| \mathcal{F}_{r_n^\delta}\right] \\
	&= e^{-(1+\frac{\Vert\theta\Vert^2}{2})r_n^\delta} \sum_{u\in N(r_n^\delta)} e^{-\theta \cdot \mathbf{X}_u(r_n^\delta)} \Pi_\mathbf{0}
	\left(\mathbf{B}_{r_n -r_n^\delta}+ \mathbf{y}+ \theta r_n^\delta \in A_n\right)
	\big|_{\mathbf{y}=\mathbf{X}_u(r_n^\delta)}.
\end{align}
Therefore,
\begin{align}
	\mu_{r_n}^\theta \left(A_n\right) - \mathbb{E} \left[ \mu_{r_n}^\theta \left(A_n\right) \big| \mathcal{F}_{r_n^\delta}\right] =: e^{-(1+\frac{\Vert\theta\Vert^2}{2})r_n^\delta} \sum_{u\in N(r_n^\delta) } e^{-\theta\cdot \mathbf{X}_u(r_n^\delta)}Y_{n,u},
\end{align}
where
\begin{align}
	Y_{n,u} & : = e^{-(1+\frac{\Vert\theta\Vert^2}{2})(r_n - r_n^\delta)}
	\sum_{v\in N(r_n ): u\leq v}
	e^{-\theta \cdot \left(\mathbf{X}_v(r_n )-\mathbf{X}_u(r_n^\delta ) \right)} 1_{\{\mathbf{X}_v(r_n ) + \theta r_n \in A_n\}}\\
	& \quad - \Pi_\mathbf{0} \left(\mathbf{B}_{r_n -r_n^\delta}+ \mathbf{y}+\theta r_n^\delta \in A_n\right) \big|_{\mathbf{y}=\mathbf{X}_u(r_n^\delta)}.
\end{align}
By the branching property, we see that, conditioned on $\mathcal{F}_{r_n^\delta}$, $\{Y_{n,u}: u\in N(r_n^\delta) \}$ is a family of centered independent random variables.  Furthermore, it holds that
\begin{align}\label{step_2}
	|Y_{n,u}|& \leq e^{-(1+\frac{\Vert\theta\Vert^2}{2})(r_n - r_n^\delta)}
	\sum_{v\in N(r_n ): u\leq v}
	e^{-\theta \cdot \left(\mathbf{X}_v(r_n )-\mathbf{X}_u(r_n^\delta ) \right)} +
	1\nonumber\\
	& = W_{r_n-r_n^\delta}(\theta ;u) +1.
\end{align}
Therefore, the second result is valid by (i) by taking $\beta =1/\kappa$ and $a_n= r_n -r_n^\delta$.
\hfill$\Box$
\bigskip

Now we treat the case $ s\in[r_n, r_{n+1})$. We will take $\delta=1/2, \beta =1/\kappa$ and $\alpha= m/2$ for $m\in \mathbb{N}$. Then the condition $\lambda \delta -\alpha >\kappa(1+\beta)$ is equivalent to $\lambda > m+2(\kappa+1).$

\begin{lemma}\label{lemma2}
	For $\mathbf{b} \in \R^d$, let $\mathbf{b}_s:= \mathbf{b}\sqrt{s}$ or $\mathbf{b}_s:= \mathbf{b}$.
	For any given $m\in \mathbb{N}$, assume that $\kappa > m +2$ and that \eqref{LlogL} holds for some $\lambda> m+2(\kappa +1)$.
	Define $k_s:=\sqrt{ r_n}$  for $s\in [r_n, r_{n+1})$.
	Then for any $\mathbf{b}\in \mathbb{R}^d$,
	\begin{align}
		s^{m/2}\left| \mu_{s}^\theta \left((-\infty, \mathbf{b}_s]\right) - \mathbb{E} \left[ \mu_{s}^\theta \left((-\infty, \mathbf{b}_s]\right) \big| \mathcal{F}_{k_s}\right]  \right| \stackrel{s\to\infty}{\longrightarrow} 0,\quad
		\mathbb{P}\mbox{-a.s.}
	\end{align}
\end{lemma}
\textbf{Proof:}
\textbf{Step 1: } In this step, we prove that almost surely,
\begin{align}\label{step_11}
	\sup_{ r_n \leq s< r_{n+1}} s^{m/2} \left| \mathbb{E} \left[\mu_{s}^\theta \left((-\infty, \mathbf{b}_s]\right) \big| \mathcal{F}_{k_s}\right] - \mathbb{E} \left[ \mu_{r_n}^\theta \left((-\infty, \mathbf{b}_{r_n}]\right) \big| \mathcal{F}_{k_s}\right]\right|\stackrel{n\to\infty}{\longrightarrow}0.
\end{align}
By the Markov property and Lemma \ref{General-many-to-one}, we see that
\begin{align}
	& \mathbb{E} \left[ \mu_{s}^\theta \left((-\infty, \mathbf{b}_s]\right) \big| \mathcal{F}_{k_s}\right] \\
	& = e^{-(1+\frac{\Vert\theta\Vert^2}{2})\sqrt{r_n}} \sum_{u\in N(\sqrt{r_n}) } e^{-\theta \cdot \mathbf{X}_u(\sqrt{r_n})}  \Pi_\mathbf{0} \left(\mathbf{B}_{s -\sqrt{r_n}} + \mathbf{y}+ \theta \sqrt{r_n} \leq \mathbf{b}_s \right)\Big|_{\mathbf{y}=\mathbf{X}_u(\sqrt{r_n}) }\\
	& = e^{-(1+\frac{\Vert\theta\Vert^2}{2})\sqrt{r_n}} \sum_{u\in N(\sqrt{r_n}) } e^{-\theta\cdot  \mathbf{X}_u(\sqrt{r_n})}  \Phi_d\Big(\frac{\mathbf{b}_s-\theta \sqrt{r_n} -\mathbf{X}_u(\sqrt{r_n})}{\sqrt{s-\sqrt{r_n}}} \Big),
\end{align}
Thus, for $s\in [r_n, r_{n+1})$, it holds that
\begin{align}
	&s^{m/2} \left| \mathbb{E} \left[\mu_{s}^\theta \left((-\infty, b_s]\right) \big| \mathcal{F}_{k_s}\right] - \mathbb{E} \left[ \mu_{r_n}^\theta \left((-\infty, b_{r_n}]\right) \big| \mathcal{F}_{k_s}\right]\right| \nonumber \\
	& \leq  r_{n+1}^{m/2} e^{-(1+\frac{\Vert\theta\Vert^2}{2})\sqrt{r_n}} \sum_{u\in N(\sqrt{r_n}) }  e^{-\theta \cdot \mathbf{X}_u(\sqrt{r_n})} \nonumber\\
	&\quad\quad  \times  \Big|\Phi_d\Big(\frac{\mathbf{b}_s-\theta \sqrt{r_n} -\mathbf{X}_u(\sqrt{r_n})}{\sqrt{s-\sqrt{r_n}}} \Big) - \Phi_d\Big(\frac{\mathbf{b}_{r_n}-\theta \sqrt{r_n} -\mathbf{X}_u(\sqrt{r_n})}{\sqrt{r_n-\sqrt{r_n}}} \Big)\Big| \nonumber \\
	&=:r_{n+1}^{m/2} e^{-(1+\frac{\Vert\theta\Vert^2}{2})\sqrt{r_n}} \sum_{u\in N(\sqrt{r_n}) }  e^{-\theta \cdot \mathbf{X}_u(\sqrt{r_n})} R(u, s).
\end{align}
Note that, on the event $\cup_{j=1}^d\{|(\mathbf{X}_u(\sqrt{r_n}))_j+\theta_j \sqrt{r_n}|> \sqrt{r_n}\}$,
we use the trivial upper-bound
$\sup_{ r_n \leq s< r_{n+1}}R(u,s) \leq 2$.
Using Lemma \ref{General-many-to-one} and the fact that $\{\mathbf{X}_\xi (t)+\theta t, \mathbb{P}^{-\theta}\}$ is a $d$-dimensional  standard Brownian motion, we have
\begin{align}
	& \sum_{n=1}^\infty r_{n+1}^{m/2} e^{-(1+\frac{\Vert\theta\Vert^2}{2})\sqrt{r_n}}  \mathbb{E}\Big(\sum_{u\in N(\sqrt{r_n}) }  e^{-\theta \cdot \mathbf{X}_u(\sqrt{r_n})} 1_{\{\cup_{j=1}^d\{|(\mathbf{X}_u(\sqrt{r_n}))_j+\theta_j \sqrt{r_n}|> \sqrt{r_n}\}\}} \sup_{ r_n \leq s< r_{n+1}} R(u, s) \Big) \\
	& \leq   2 \sum_{n=1}^\infty r_{n+1}^{m/2} \Pi_\mathbf{0} \left(\cup_{j=1}^d\left\{|(\mathbf{B}_1)_j| > r_n^{1/4}\right\}\right)\leq 2d \Pi_\mathbf{0} (e^{|(\mathbf{B}_1)_1|}) \sum_{n=1}^\infty r_{n+1}^{m/2} e^{-r_n^{1/4}}<\infty,
\end{align}
which implies that $\P$-almost surely,
\begin{align}\label{step_9}
	r_{n+1}^{m/2}  e^{-(1+\frac{\Vert\theta\Vert^2}{2})\sqrt{r_n}} \sum_{u\in N(\sqrt{r_n}) }  e^{-\theta\cdot \mathbf{X}_u(\sqrt{r_n})} 1_{\{\cup_{j=1}^d\{|(\mathbf{X}_u(\sqrt{r_n}))_j+\theta_j \sqrt{r_n}|> \sqrt{r_n}\}\}}R(u, s) \stackrel{s\to\infty}{\longrightarrow} 0.
\end{align}
On the other hand, on the event $\cap_{j=1}^d\{|(\mathbf{X}_u(\sqrt{r_n}))_j+\theta_j \sqrt{r_n}|\leq \sqrt{r_n}\}$,
in the case $\mathbf{b}_s= \mathbf{b}\sqrt{s}$, using the trivial inequality
\[
\left|\Phi_d(\mathbf{a})- \Phi_d(\mathbf{b})\right| \leq \sum_{j=1}^d \left|\Phi(a_j)- \Phi(b_j)\right|\leq\frac{1}{\sqrt{2\pi}}\sum_{j=1}^d|a_j-b_j|,
\]
we get that, uniformly for $s\in [r_n r_{n+1})$,
\begin{align}
	& R(u,s)\leq \sum_{j=1}^d \Bigg(	\frac{|b_j|}
	{\sqrt{2\pi}} \left| \frac{\sqrt{s}}{\sqrt{s-\sqrt{r_n}}} -\frac{\sqrt{r_n}}{\sqrt{r_n-\sqrt{r_n}}} \right| \\
	&\quad\quad + \frac{|(\mathbf{X}_u(\sqrt{r_n}))_j+\theta_j \sqrt{r_n}| }{\sqrt{2\pi}} \left|\frac{1}{\sqrt{s-\sqrt{r_n}}} - \frac{1}{\sqrt{r_n -\sqrt{r_n}}}\right| \Bigg)\nonumber\\
	& \lesssim  \frac{\left|\sqrt{s(r_n -\sqrt{r_n}) }- \sqrt{r_n(s-\sqrt{r_n})} \right|}{\sqrt{s-\sqrt{r_n}} \sqrt{r_n -\sqrt{r_n}}} + \sqrt{r_n}\frac{\left|\sqrt{s-\sqrt{r_n}} -\sqrt{r_n -\sqrt{r_n}}\right|}{ \sqrt{s-\sqrt{r_n}} \sqrt{r_n -\sqrt{r_n}}}\\
	& \lesssim \frac{1}{r_n} \frac{\left|s(r_n -\sqrt{r_n}) - r_n(s-\sqrt{r_n})\right|}{r_n} + \frac{\sqrt{r_n}}{r_n^{3/2}}(r_{n+1}-r_n)\lesssim \frac{1}{r_n}(r_{n+1}-r_n).
\end{align}
In the case $\mathbf{b}_s=\mathbf{b}$, we have
\begin{align}
	&	R(u,s)\leq \sum_{j=1}^d\Bigg( \frac{|b_j|}
	{\sqrt{2\pi}} \left| \frac{1}{\sqrt{s-\sqrt{r_n}}} -\frac{1}{\sqrt{r_n-\sqrt{r_n}}} \right| \\
	&\quad\quad+ \frac{|(\mathbf{X}_u(\sqrt{r_n}))_j+\theta_j \sqrt{r_n}| }{\sqrt{2\pi}} \left|\frac{1}{\sqrt{s-\sqrt{r_n}}} - \frac{1}{\sqrt{r_n -\sqrt{r_n}}}\right|\Bigg)\nonumber\\
	& \lesssim \frac{\left|s-r_n \right|}{\sqrt{s-\sqrt{r_n}} \sqrt{r_n -\sqrt{r_n}} \left(\sqrt{s-\sqrt{r_n}}+\sqrt{r_n-\sqrt{r_n}} \right)} + \sqrt{r_n}\frac{\left|\sqrt{s-\sqrt{r_n}} -\sqrt{r_n -\sqrt{r_n}}\right|}{ \sqrt{s-\sqrt{r_n}} \sqrt{r_n -\sqrt{r_n}}}\\
	& \lesssim  \frac{1}{r_n^{3/2}}(r_{n+1}-r_n)+ \frac{\sqrt{r_n}}{r_n^{3/2}}(r_{n+1}-r_n)\lesssim \frac{1}{r_n}(r_{n+1}-r_n).
\end{align}
Thus in both cases, we have that
\begin{align}\label{step_10}
	& r_{n+1}^{m/2}  e^{-(1+\frac{\Vert\theta\Vert^2}{2})\sqrt{r_n}} \sum_{u\in N(\sqrt{r_n}) }  e^{-\theta \cdot \mathbf{X}_u(\sqrt{r_n})} 1_{\{\cap_{j=1}^d\{|(\mathbf{X}_u(\sqrt{r_n}))_j+\theta_j \sqrt{r_n}|\leq \sqrt{r_n}\}\}}
	\sup_{ r_n \leq s< r_{n+1}} R(u, s) \nonumber\\
	& \lesssim r_{n+1}^{m/2} \frac{1}{r_n}(r_{n+1}-r_n)  W_{\sqrt{r_n}}(\theta).
\end{align}
We claim that right-hand side of \eqref{step_10} goes to 0 almost surely as $s\to\infty$.
In fact,
\[
r_{n+1}^{m/2} \frac{1}{r_n}(r_{n+1}-r_n) \lesssim
r_{n}^{(-2+m)/2} (r_{n+1}-r_n)
=n^{(-2+m)/(2\kappa)}\left((n+1)^{1/\kappa} - n^{1/\kappa}\right).
\]
By the mean value theorem, the right-hand side above is equal to
\[
n^{(-2+m)/(2\kappa)} \xi^{-1+\frac{1}{\kappa}}\lesssim n^{(m-2\kappa)/(2\kappa)}.
\]
Since $\kappa >m+2$, the claim is valid.
Combining this with \eqref{step_9} and \eqref{step_10}, we get \eqref{step_11}.

\textbf{Step 2: } In this step, we prove that
\begin{align}\label{step_12}
	\sup_{ r_n \leq s< r_{n+1}}	s^{m/2}\left| \mu_{s}^\theta \left((-\infty, \mathbf{b}_s]\right) - \mu_{r_n}^\theta \left((-\infty, \mathbf{b}_{r_n}]\right) \right| \stackrel{n\to\infty}{\longrightarrow} 0,\quad
	\mathbb{P}\mbox{-a.s.}
\end{align}
Once we get \eqref{step_12}, we can
combine
\eqref{step_11} and Lemma \ref{lemma1} (ii)
(with $A_n= (-\infty, \mathbf{b}_{r_n}]$ and $\delta=1/2$)
to get the assertion of the lemma.

To prove \eqref{step_12}, we first prove that
\begin{align}\label{step_18}
	\liminf_{n\to\infty}  \inf_{r_n \leq s< r_{n+1}} s^{m/2} \left( \mu_{s}^\theta \left((-\infty, \mathbf{b}_s]\right)  -\mu_{r_n}^\theta \left((-\infty, \mathbf{b}_{r_n}]\right) \right)\geq 0,\quad\mathbb{P}\mbox{-a.s.}
\end{align}
Define $\varepsilon_n:= \sqrt{r_{n+1}- r_n}$.
For $u\in N(r_n)$, let $G_u$ be the event that
$u$ does not split before $r_{n+1}$ and that
$\max_{s\in (r_n, r_{n+1})}\Vert {\mathbf X}_u(s)-{\mathbf X}_u(r_n)\Vert \le \sqrt{r_n}\varepsilon_n$.
Then
\[
\mathbb{P}\left(G_u\big| \mathcal{F}_{r_n}\right) = e^{-(r_{n+1}-r_n)}\Pi_\mathbf{0} \Big(\max_{r\leq r_{n+1}-r_n}\Vert\mathbf{B}_r\Vert \leq \sqrt{r_n} \varepsilon_n\Big) =e^{-(r_{n+1}-r_n)}\Pi_\mathbf{0} \Big(\max_{r\leq 1}\Vert\mathbf{B}_r\Vert \leq \sqrt{r_n}\Big) .
\]
Recalling that $\mathbf{1}:=(1,...,1)$, it holds that
\begin{align}\label{Lower-bound}
	& \mu_{s}^\theta \left((-\infty, {\mathbf b}_s]\right)
	= e^{-(1+\frac{\Vert\theta\Vert^2}{2})s}\sum_{u\in N(r_n)} e^{-\theta\cdot \mathbf{X}_u(r_n)}
	\sum_{v\in N(s): u\leq v}
	e^{-\theta \cdot \left(\mathbf{X}_v(s)-\mathbf{X}_u(r_n) \right)} 1_{\{\mathbf{X}_v(s)+\theta s \leq \mathbf{b}_s\}}\nonumber\\
	& \geq e^{-(1+\frac{\Vert\theta\Vert^2}{2})r_{n+1}} e^{-\Vert\theta\Vert \sqrt{r_n}\varepsilon_n}\sum_{u\in N(r_n)} e^{-\theta \cdot \mathbf{X}_u(r_n)} 1_{\{\mathbf{X}_u(r_n) +\theta r_n \leq \mathbf{b}_{r_n}-\varepsilon_n\sqrt{r_n} \mathbf{1} - \Vert \theta \Vert (r_{n+1}-r_n) \mathbf{1} \} } 1_{G_u}\nonumber \\
	& = e^{-(1+\frac{\Vert\theta\Vert^2}{2})r_{n+1}-\Vert\theta\Vert \sqrt{r_n}\varepsilon_n}\sum_{u\in N(r_n)} e^{-\theta \cdot \mathbf{X}_u(r_n)} \nonumber\\
	& \quad \quad \times 1_{\{\mathbf{X}_u(r_n) +\theta r_n \leq \mathbf{b}_{r_n}-\varepsilon_n\sqrt{r_n}\mathbf{1}-\Vert \theta \Vert (r_{n+1}-r_n) \mathbf{1} \} } \left(1_{G_u}- \mathbb{P}\left(G_u\big| \mathcal{F}_{r_n}\right) \right)\nonumber \\
	&\quad + e^{-(1+\frac{\Vert\theta\Vert^2}{2})r_{n+1}-\Vert\theta\Vert \sqrt{r_n}\varepsilon_n}\sum_{u\in N(r_n)} e^{-\theta \cdot \mathbf{X}_u(r_n)} \nonumber\\
	& \quad \quad \times 1_{\{\mathbf{X}_u(r_n) +\theta r_n \leq \mathbf{b}_{r_n}-\varepsilon_n\sqrt{r_n}\mathbf{1}-\Vert \theta \Vert (r_{n+1}-r_n) \mathbf{1} \} }  \mathbb{P}\left(G_u\big| \mathcal{F}_{r_n}\right) =: I + II.
\end{align}
For $I$, we will apply Lemma \ref{lemma1} (i) with
$\alpha=m/2$, $\delta =1, a_n=0, \beta =0$ and
$$
Y_{n,u}:=
1_{\{\mathbf{X}_u(r_n) +\theta r_n \leq \mathbf{b}_{r_n}-\varepsilon_n\sqrt{r_n}\mathbf{1}-\Vert \theta \Vert (r_{n+1}-r_n) \mathbf{1} \} }
\left(1_{G_u}- \mathbb{P}\left(G_u\big| \mathcal{F}_{r_n}\right) \right).
$$
It is easy to see that $|Y_{n,u}|\leq 2$, $r_{n+1}-r_n \to 0$ and
$\sqrt{r_n}\varepsilon_n \lesssim \sqrt{n^{(2-\kappa)/\kappa}}\to 0$.
Since $\lambda>m+2(\kappa+1)$, we have
\begin{align}\label{step_16}
	\sup_{ r_n \leq s< r_{n+1}} s^{m/2} |I| \stackrel{n\to\infty}{\longrightarrow} 0,\quad
	\mathbb{P}\mbox{-a.s.}
\end{align}
If we can prove that
\begin{align}\label{step_15}
	\sup_{ r_n \leq s< r_{n+1}}s^{m/2} \left|II -  \mu_{r_n}^\theta \left((-\infty, \mathbf{b}_{r_n}]\right) \right|\stackrel{n\to\infty}{\longrightarrow} 0,\quad
	\P\mbox{-a.s.},
\end{align}
then \eqref{step_18} will follow from \eqref{Lower-bound},  \eqref{step_16} and  \eqref{step_15}.
Now we prove  \eqref{step_15}. Since $\kappa >m+2$, we have
$r_n^{m/2}(r_{n+1}-r_n)\lesssim n^{-1+(m+2)/(2\kappa) }\to 0$.
Thus,
\begin{align}
	& s^{m/2} \left|1-\mathbb{P}\left(G_u\big| \mathcal{F}_{r_n}\right) \right|\leq r_{n+1}^{m/2} (1- e^{-(r_{n+1}-r_n)}) + r_{n+1}^{m/2} \Pi_\mathbf{0} \Big(\max_{r\leq 1}\Vert \mathbf{B}_r\Vert > \sqrt{r_n}\Big) \\
	& \lesssim r_{n+1}^{m/2} (r_{n+1} -r_n) + r_{n+1}^{m/2}e^{-\sqrt{r_n}} \to 0.
\end{align}
Hence,
\begin{align}\label{step_13}
	& s^{m/2} \left|e^{(1+\frac{\Vert\theta\Vert^2}{2})(r_{n+1}-r_n)} e^{\Vert\theta\Vert \sqrt{r_n}\varepsilon_n}II -  \mu_{r_n}^\theta \left((-\infty, \mathbf{b}_{r_n}-\varepsilon_n  \sqrt{r_n} \mathbf{1}
	-\Vert \theta \Vert (r_{n+1}-r_n) \mathbf{1}]\right) \right| \nonumber\\
	& \lesssim W_{r_n}(\theta ) \left(r_{n+1}^{m/2} (r_{n+1} -r_n) + r_{n+1}^{m/2} e^{-\sqrt{r_n}} \right) \stackrel{s\to\infty}{\longrightarrow} 0,\quad
	\mathbb{P}\mbox{-a.s.}
\end{align}
Note that $0\leq e^{(1+\frac{\Vert\theta\Vert^2}{2})(r_{n+1}-r_n)}e^{\Vert\theta\Vert \sqrt{r_n}\varepsilon_n}II \leq W_{r_n}(\theta)$, $e^{(1+\frac{\Vert\theta\Vert^2}{2})(r_{n+1}-r_n)} = 1+ O(r_{n+1}-r_n)= 1+ o(r_n^{-m/2})$, and that  $e^{\Vert\theta\Vert \sqrt{r_n} \varepsilon_n} = 1+ O(\sqrt{r_n} \sqrt{r_{n+1}-r_n}) = 1+O(n^{1/\kappa - 1/2}) = 1+o(r_n^{-m/2})$ by the assumption that $\kappa > m+2$. Therefore, \eqref{step_13} implies that
\begin{align}\label{step_14}
	s^{m/2} \left|II -  \mu_{r_n}^\theta \left((-\infty, \mathbf{b}_{r_n}-\varepsilon_n \sqrt{r_n}\mathbf{1}
	-\Vert \theta \Vert (r_{n+1}-r_n) \mathbf{1}]\right)
	\right|\stackrel{s\to\infty}{\longrightarrow} 0,\quad
	\mathbb{P}\mbox{-a.s.}
\end{align}
Now we put
$A_n=  (-\infty, \mathbf{b}_{r_n}] \setminus (-\infty, \mathbf{b}_{r_n}-\varepsilon_n \sqrt{r_n} \mathbf{1} - \Vert \theta \Vert (r_{n+1}-r_n) \mathbf{1}]\subset \cup_{j=1}^d  C_{n,j}$ where
$C_{n,j}:= \left\{ \mathbf{x}=(x_1,...,x_d): x_j \in ((\mathbf{b}_{r_n})_j -\varepsilon_n \sqrt{r_n} - \Vert \theta \Vert(r_{n+1}-r_n),(\mathbf{b}_{r_n})_j ] \right\}$.
Then by Lemma \ref{General-many-to-one} and the inequality $\Pi_\mathbf{0}\left(\mathbf{B}_t + \mathbf{y}\in C_{n,j}\right)\leq
\frac{\varepsilon_n \sqrt{r_n}+ \Vert\theta \Vert (r_{n+1}-r_n)}{\sqrt{2\pi t}}$, we obtain that
\begin{align}
	& r_n^{m/2} \mathbb{E} \left[ \mu_{r_n}^\theta \left(A_n\right) \big| \mathcal{F}_{\sqrt{r_n}}\right] \\
	& =  r_n^{m/2}e^{-(1+\frac{\Vert \theta\Vert^2}{2})\sqrt{r_n}}   \sum_{u\in N(\sqrt{r_n})} e^{-\theta\cdot \mathbf{X}_u(\sqrt{r_n})} \P^{-\theta}\left(\mathbf{X}_\xi(r_n - \sqrt{r_n}) +\mathbf{y} +\theta r_n  \in A_n \right)|_{\mathbf{y}=\mathbf{X}_u(\sqrt{r_n}) }\\
	& \leq \sum_{j=1}^d  r_n^{m/2}e^{-(1+\frac{\Vert \theta\Vert^2}{2})\sqrt{r_n}}   \sum_{u\in N(\sqrt{r_n})} e^{-\theta\cdot \mathbf{X}_u(\sqrt{r_n})} \Pi_\mathbf{0}\left(\mathbf{B}_{r_n - \sqrt{r_n}} +\mathbf{y} +\theta \sqrt{r_n}  \in C_{n,j} \right)|_{\mathbf{y}=\mathbf{X}_u(\sqrt{r_n}) }\\
	& \leq d  W_{\sqrt{r_n}}(\theta) \frac{r_n^{m/2}
		\left( \varepsilon_n \sqrt{r_n}+ \Vert\theta \Vert (r_{n+1}-r_n)\right)
	}{\sqrt{2\pi (r_n - \sqrt{r_n})}} \stackrel{s\to\infty}{\longrightarrow} 0.
\end{align}
Here the last assertion about the limit being 0 follows from the following argument:
\[
\frac{r_n^{m/2}
	\left( \varepsilon_n \sqrt{r_n}+ \Vert\theta \Vert (r_{n+1}-r_n)\right)
}{\sqrt{2\pi (r_n - \sqrt{r_n})}} \lesssim r_n^{m/2} \varepsilon_n = n^{m /(2\kappa)}\sqrt{(n+1)^{1/\kappa} - n^{1/\kappa}}\lesssim
n^{ (m+1-\kappa)/(2\kappa)}\to 0.
\]
Using Lemma \ref{lemma1} (ii), we immediately get that $r_n^{m/2} \mu_{r_n}^\theta (A_n) \to 0$, $\mathbb{P}$-almost surely. Then by \eqref{step_14},
we conclude that \eqref{step_15} holds.

Applying similar arguments  for the interval $(\mathbf{b}_s, +\infty)$, we can also get
\begin{align}\label{step_19}
	\liminf_{n\to\infty}  \inf_{r_n \leq s< r_{n+1}} s^{m/2} \left( \mu_{s}^\theta \left( (\mathbf{b}_s,  +\infty )\right)  -\mu_{r_n}^\theta \left((\mathbf{b}_{r_n}, +\infty )\right) \right)\geq 0,\quad \P\mbox{-a.s.}
\end{align}
Using  Proposition \ref{Convergence-rate-Martingale}  with
$\mathbf{k}=\mathbf{0}$ and $\eta=2(\kappa+1)$,  and the assumption $\lambda >m+2(\kappa+1)$, we get
\begin{align}\label{step_17}
	\lim_{n\to\infty} \sup_{ r_n \leq s< r_{n+1}} s^{m/2} \left| \mu_{s}^\theta \left( \mathbb{R}^d\right)  -\mu_{r_n}^\theta \left(\mathbb{R}^d\right) \right| = \lim_{n\to\infty} \sup_{ r_n \leq s< r_{n+1}} s^{m/2} \left|W_s(\theta) -W_{r_n}(\theta) \right| = 0.
\end{align}
Now we prove \eqref{step_12} follows from \eqref{step_18}, \eqref{step_19} and \eqref{step_17}.
Indeed,  for any $\varepsilon>0$,  \eqref{step_18}, \eqref{step_19} and \eqref{step_17} imply
that one can find a random time $N$ such that for all $n>N$ and $r_n \leq s < r_{n+1}$,
\begin{align}
	& s^{m/2} \left( \mu_{s}^{m/2} \left((-\infty, \mathbf{b}_s]\right)  -\mu_{r_n}^\theta \left((-\infty, \mathbf{b}_{r_n}]\right) \right) > -\varepsilon,\\
	& s^{m/2} \left( \mu_{s}^\theta \left( (\mathbf{b}_s, +\infty )\right)  -\mu_{r_n}^\theta \left((\mathbf{b}_{r_n}, +\infty )\right) \right) > -\varepsilon\quad \mbox{and}\quad s^{m/2} \left| \mu_{s}^\theta \left( \mathbb{R}^d\right)  -\mu_{r_n}^\theta \left(\mathbb{R}^d\right) \right| < \varepsilon.
\end{align}
Thus,
\begin{align}
	& s^{m/2} \left( \mu_{s}^\theta \left((-\infty, \mathbf{b}_s]\right)  -\mu_{r_n}^\theta \left((-\infty, \mathbf{b}_{r_n}]\right) \right) \\
	&= s^{m/2} \left( \mu_{s}^\theta \left( \mathbb{R}^d\right)  -\mu_{r_n}^\theta \left(\mathbb{R}^d\right) \right) - s^{m/2} \left( \mu_{s}^\theta \left( (\mathbf{b}_s, +\infty )\right)  -\mu_{r_n}^\theta \left((\mathbf{b}_{r_n}, +\infty )\right) \right) <2\varepsilon.
\end{align}
Hence we have that when $n>N$ and $r_n\leq s < r_{n+1}$,
\[
s^{m/2} \left| \mu_{s}^\theta \left((-\infty, \mathbf{b}_s]\right)  -\mu_{r_n}^\theta \left((-\infty, \mathbf{b}_{r_n}]\right) \right|< 2\varepsilon,
\]
which implies  \eqref{step_12}.
\hfill$\Box$
\bigskip

For any given $m\in \mathbb{N}$, we will take $\kappa:= m+3$ in the remainder of this section. It follows from Lemma \ref{lemma2} that, if \eqref{LlogL} holds for some $\lambda > m+ 2(\kappa +1)= 3m+8$, then $\mathbb{P}$-almost surely for all $s\in [r_n, r_{n+1})$ and $\mathbf{b}_s= \mathbf{b}\sqrt{s}$,
\begin{align}\label{step_20}
	& \mu_{s}^\theta \left((-\infty, \mathbf{b}\sqrt{s}]\right) = \mathbb{E} \left[ \mu_{s}^\theta \left((-\infty, \mathbf{b}\sqrt{s}]\right) \big| \mathcal{F}_{k_s}\right]  + o(s^{-m/2})\nonumber\\
	& =  e^{-(1+\frac{\Vert\theta\Vert^2}{2})\sqrt{r_n}} \sum_{u\in N(\sqrt{r_n}) } e^{-\theta\cdot \mathbf{X}_u(\sqrt{r_n})}  \Phi_d\Big(\frac{\mathbf{b}\sqrt{s}-\theta \sqrt{r_n} -\mathbf{X}_u(\sqrt{r_n})}{\sqrt{s-\sqrt{r_n}}} \Big)\nonumber\\
&\quad + o(s^{-m/2}).
\end{align}
Note that, for any $\mathbf{a}<\mathbf{b}$,
$(\mathbf{a}, \mathbf{b}]= \prod_{j=1}^d (a_j, b_j]$
can be expressed
in terms $\prod_{j=1}^d E_j$ where
$E_j\in \{(-\infty, a_j],  (-\infty, b_j] \}$ using a finite number of set theoretic operations.
Thus, applying Lemma \ref{lemma2} to $\prod_{j=1}^d E_j$, we get that, if  \eqref{LlogL} holds for some $\lambda> 3(m+d)+8= 3m+3d+8$, then
\begin{align}\label{step_24}
	& \mu_s^\theta \left( (\mathbf{a},\mathbf{b}]\right)
	=   o(s^{-(m+d)/2})+ e^{-(1+\frac{\Vert\theta\Vert^2}{2})\sqrt{r_n}} \sum_{u\in N(\sqrt{r_n}) } e^{-\theta \cdot \mathbf{X}_u(\sqrt{r_n})}\nonumber\\
	& \quad\quad \quad \times
	\prod_{j=1}^d \frac{1}{\sqrt{ s-\sqrt{r_n}}} \int_{a_j}^{b_j}  \phi\Big( \frac{z_j-\theta_j \sqrt{r_n}- (\mathbf{X}_u(\sqrt{r_n}))_j}{\sqrt{s-r_n}}\Big)\mathrm{d}z_j.
\end{align}

\bigskip

\textbf{Proof of Theorem \ref{thm1}:}
Let $m\in \N$ and assume \eqref{LlogL} holds for some $\lambda> \max\left\{ 3m+8, d(3m+5)\right\}$.
Recall that
$r_n=n^{1/\kappa}$
and  $\kappa=m+3$. Put $K:= m/\kappa +3$.
Combining Lemma \ref{General-many-to-one}, $\sup_{\mathbf{z}\in \R^d}\Phi_d(\mathbf{z})= 1$ and the fact that $\{(\mathbf{X}_\xi(t)+\theta t)_{t\geq 0},\mathbb{P}^{-\theta}\} $ is a $d$-dimensional standard Brownian motion, we get that
\begin{align}\label{step_26}
	& \sum_{n=2}^\infty
	r_n^{m/2} \mathbb{E} \Big(	e^{-(1+\frac{\Vert\theta\Vert^2}{2})\sqrt{r_n}} \sum_{u\in N(\sqrt{r_n}) } e^{-\theta \cdot \mathbf{X}_u(\sqrt{r_n})} \nonumber\\
	&\quad \quad \times \sup_{ r_n \leq s< r_{n+1}} \Phi_d\Big(\frac{\mathbf{b}\sqrt{s}-\theta \sqrt{r_n} -\mathbf{X}_u(\sqrt{r_n})}{\sqrt{s-\sqrt{r_n}}} \Big) 1_{\left\{ \cup_{j=1}^d\left\{ \left|(\mathbf{X}_u(\sqrt{r_n}) )_j+\theta_j \sqrt{r_n}\right| > \sqrt{K\sqrt{r_n} \log n} \right\} \right\}} \Big)\nonumber\\
	& \leq   \sum_{n=2}^\infty
	n^{m/(2\kappa)} \sum_{j=1}^d \Pi_\mathbf{0} \left(|(\mathbf{B}_{\sqrt{r_n}})_j| > \sqrt{K \sqrt{r_n} \log n}\right) \nonumber\\
	& = d \sum_{n=2}^\infty
	n^{m/(2\kappa)} \Pi_\mathbf{0} \left(|(\mathbf{B}_1)_1| > \sqrt{K \log n}\right)\lesssim \sum_{n=1}^\infty n^{m/(2\kappa)} n^{-K/2}<\infty,
\end{align}
where in the last inequality we used  the fact that $\Pi_\mathbf{0}(|(\mathbf{B}_1)_1| >x)\lesssim e^{-x^2/2}$.
Therefore, $\P$-almost surely,
\begin{align}\label{limit0}
	&r_n^{m/2} 	e^{-(1+\frac{\Vert\theta\Vert^2}{2})\sqrt{r_n}} \sum_{u\in N(\sqrt{r_n}) } e^{-\theta \cdot \mathbf{X}_u(\sqrt{r_n})}1_{\left\{ \cup_{j=1}^d\left\{ \left|(\mathbf{X}_u(\sqrt{r_n}) )_j+\theta_j \sqrt{r_n}\right| > \sqrt{K\sqrt{r_n} \log n} \right\} \right\}}  \nonumber \\
	&\quad \quad\times \sup_{ r_n \leq s< r_{n+1}} \Phi_d\Big(\frac{\mathbf{b}\sqrt{s}-\theta \sqrt{r_n} -\mathbf{X}_u(\sqrt{r_n})}{\sqrt{s-\sqrt{r_n}}} \Big) \stackrel{n\to\infty}{\longrightarrow} 0.
\end{align}
Since $\lambda>3m+8$,
by \eqref{step_20} and  \eqref{limit0},  for any $\theta\in \R^d$ with $\Vert\theta\Vert< \sqrt{2},
\mathbf{b}\in \mathbb{R}^d$ and $s\in [r_n, r_{n+1})$,
\begin{align}\label{step_21}
	& \mu_{s}^\theta \left((-\infty, \mathbf{b}\sqrt{s}]\right) =
	o(s^{-m/2}) +e^{-(1+\frac{\Vert\theta\Vert^2}{2})\sqrt{r_n}} \sum_{u\in N(\sqrt{r_n}) } e^{-\theta \cdot \mathbf{X}_u(\sqrt{r_n})} \nonumber\\
	& \quad \quad \times \Phi_d\Big(\frac{\mathbf{b}\sqrt{s}-\theta \sqrt{r_n} -\mathbf{X}_u(\sqrt{r_n})}{\sqrt{s-\sqrt{r_n}}} \Big) 1_{\left\{ \cap_{j=1}^d \left\{ \left|(\mathbf{X}_u(\sqrt{r_n}) )_j+\theta_j \sqrt{r_n}\right| \leq  \sqrt{K\sqrt{r_n} \log n}\right\} \right\}}.
\end{align}
Put $J:= 6m+10$. Then  $J> 2m+ K\kappa = 3m+3\kappa = 6m+9 $.
By Lemma \ref{lemma3}, we get that for any  $\theta\in \R^d$ with $\Vert\theta\Vert< \sqrt{2}, \mathbf{b}\in \mathbb{R}^d$ and $s\in[r_n, r_{n+1})$,
\begin{align}\label{step_21'}
	& \mu_{s}^\theta \left((-\infty, \mathbf{b}\sqrt{s}]\right)\nonumber\\
	& = o(s^{-m/2}) +e^{-(1+\frac{\Vert\theta\Vert^2}{2})\sqrt{r_n}} \sum_{u\in N(\sqrt{r_n}) } e^{-\theta \cdot \mathbf{X}_u(\sqrt{r_n})} 1_{\left\{ \cap_{j=1}^d \left\{ \left|(\mathbf{X}_u(\sqrt{r_n}) )_j+\theta_j \sqrt{r_n}\right| \leq  \sqrt{K\sqrt{r_n} \log n}\right\} \right\}} \nonumber\\
	&\quad \times  \prod_{j=1}^d\Big( \Phi (b_j)- \phi(b_j) \sum_{k=1}^J \frac{1}{k!} \frac{1}{s^{k/2}} H_{k-1}(b_j) \left(\left(\sqrt{r_n}\right)^{k/2} H_k\left(\frac{(\mathbf{X}_u(\sqrt{r_n}))_j+\theta_j \sqrt{r_n}}{r_n^{1/4}}\right)\right)+ \varepsilon_{m,u, s,j}  \Big)\nonumber\\
	&=  o(s^{-m/2}) +e^{-(1+\frac{\Vert\theta\Vert^2}{2})\sqrt{r_n}} \sum_{u\in N(\sqrt{r_n}) } e^{-\theta \cdot \mathbf{X}_u(\sqrt{r_n})} 1_{\left\{ \cap_{j=1}^d \left\{ \left|(\mathbf{X}_u(\sqrt{r_n}) )_j+\theta_j \sqrt{r_n}\right| \leq  \sqrt{K\sqrt{r_n} \log n}\right\} \right\}} \nonumber\\
	&\quad \times  \prod_{j=1}^d\Big( \Phi (b_j)- \phi(b_j) \sum_{k=1}^J \frac{1}{k!} \frac{1}{s^{k/2}} H_{k-1}(b_j) \Big(\left(\sqrt{r_n}\right)^{k/2} H_k\Big(\frac{(\mathbf{X}_u(\sqrt{r_n}))_j+\theta_j \sqrt{r_n}}{r_n^{1/4}}\Big)\Big)  \Big),
\end{align}
where $\varepsilon_{m,u, s,j}=\varepsilon_{m,y, b,s}|_{y=\theta_j \sqrt{r_n}+(\mathbf{X}_u(\sqrt{r_n}))_j, \ b=b_j}$.
To justify the last equality, we first apply Lemma \ref{lemma3} to get that, for each $u\in N(\sqrt{r_n})$, as $s\to\infty$, $\P$-almost surely,
\begin{align}
	&s^{m/2}e^{-(1+\frac{\theta^2}{2})\sqrt{r_n}} \sum_{u\in N(\sqrt{r_n}) } e^{-\theta X_u(\sqrt{r_n})} 1_{\left\{ \cap_{j=1}^d \left\{ \left|(\mathbf{X}_u(\sqrt{r_n}) )_j+\theta_j \sqrt{r_n}\right| \leq  \sqrt{K\sqrt{r_n} \log n}\right\} \right\}} \left| \varepsilon_{m,u, s,j}\right|\\
	& \leq  W_{\sqrt{r_n}}(\theta) \times
	s^{m/2}\sup_{j\leq d}\sup_{u\in N(\sqrt{r_n})} |\varepsilon_{m,u, s,j}|   1_{\left\{ \cap_{j=1}^d \left\{ \left|(\mathbf{X}_u(\sqrt{r_n}) )_j+\theta_j \sqrt{r_n}\right| \leq  \sqrt{K\sqrt{r_n} \log n}\right\} \right\}}\to 0.
\end{align}
Then note  that
by \eqref{Upper-H_k-2} and $|H_k(x)|\lesssim |x|^k +1$,
on the set $ \cap_{j=1}^d \left\{ \left|(\mathbf{X}_u(\sqrt{r_n}) )_j+\theta_j \sqrt{r_n}\right| \leq  \sqrt{K\sqrt{r_n} \log n}\right\}$,
\begin{align}
	&|Q_{u,j}|:= \Big|\Phi (b_j)- \phi(b_j) \sum_{k=1}^J \frac{1}{k!} \frac{1}{s^{k/2}} H_{k-1}(b_j) \Big(\left(\sqrt{r_n}\right)^{k/2} H_k\Big(\frac{(\mathbf{X}_u(\sqrt{r_n}))_j+\theta_j \sqrt{r_n}}{r_n^{1/4}}\Big)\Big)  \Big|
	\nonumber\\ & \lesssim  1+ \sum_{k=1}^J \frac{1}{r_n^{k/2}} r_n^{k/4} \Big(1+ \Big| \frac{(\mathbf{X}_u(\sqrt{r_n}))_j+\theta_j \sqrt{r_n}}{r_n^{1/4}}\Big|^k\Big)\lesssim 1+ \frac{(\log n)^{J/2}}{r_n^{1/4}}\lesssim 1.
\end{align}
Combining the two displays above, we get that,
on the set $ \cap_{j=1}^d \left\{ \left|(\mathbf{X}_u(\sqrt{r_n}) )_j+\theta_j \sqrt{r_n}\right| \leq  \sqrt{K\sqrt{r_n} \log n}\right\}$,
\begin{align}
	&\Big| \Phi_d\Big(\frac{\mathbf{b}\sqrt{s}-\theta \sqrt{r_n} -\mathbf{X}_u(\sqrt{r_n})}{\sqrt{s-\sqrt{r_n}}} \Big) - \prod_{j=1}^d Q_{u,j}\Big|\leq \sum_{j=1}^d \prod_{\ell\neq j} \left| Q_{v, \ell }\right| |\varepsilon_{m,u,s,j}|\\
	&
	\lesssim d
	\sup_{j\leq d}\sup_{u\in N(\sqrt{r_n})} |\varepsilon_{m,u, s,j}|   1_{\left\{ \cap_{j=1}^d \left\{ \left|(\mathbf{X}_u(\sqrt{r_n}) )_j+\theta_j \sqrt{r_n}\right| \leq  \sqrt{K\sqrt{r_n} \log n}\right\} \right\}},
\end{align}
which implies \eqref{step_21'}.

Let $\varepsilon\in (0,1)$ be small enough so that $K(1-\varepsilon)\geq  m/\kappa +2$.
For any $\mathbf{k}\in \N^d$ with $1\leq |\mathbf{k}|\leq J$,
using the inequality $|H_k(x)|\lesssim 1+|x|^k$ first and then Lemma \ref{General-many-to-one}, we get
\begin{align}\label{step_28}
	&  \sum_{n=2}^\infty
	r_n^{m/2} \mathbb{E} \sup_{ r_n \leq s< r_{n+1}}\Big|e^{-(1+\frac{\Vert\theta\Vert^2}{2})\sqrt{r_n}} \sum_{u\in N(\sqrt{r_n}) } e^{-\theta \cdot \mathbf{X}_u(\sqrt{r_n})} \nonumber\\
	& \quad \times  1_{\left\{ \cup_{j=1}^d\left\{ \left|(\mathbf{X}_u(\sqrt{r_n}))_j +\theta_j \sqrt{r_n}\right| > \sqrt{K\sqrt{r_n} \log n}  \right\}\right\}}
	\prod_{j=1}^d \frac{\left(\sqrt{r_n}\right)^{k_j/2} }{s^{k_j/2}}   H_{k_j}\Big(\frac{(\mathbf{X}_u(\sqrt{r_n}))_j+\theta_j \sqrt{r_n}}{r_n^{1/4}}\Big)  \Big|\nonumber\\
	& \lesssim  \sum_{n=2}^\infty
	r_n^{m/2} \Pi_\mathbf{0} \Big( \sum_{j=1}^d 1_{\left\{  \left|(\mathbf{B}_{\sqrt{r_n}})_j\right| > \sqrt{K\sqrt{r_n} \log n} \right\}} \prod_{\ell=1}^d\frac{1}{r_n^{k_\ell/4}}   \Big(  1+ \left|\frac{(\mathbf{B}_{\sqrt{r_n}})_\ell}{\sqrt{\sqrt{r_n}}}\right|^{k_\ell} \Big)\Big)\nonumber\\
	& =\sum_{j=1}^d  \sum_{n=2}^\infty
	n^{ (2m- |\mathbf{k}|)/(4\kappa)}
	\Pi_\mathbf{0 }\Big(  1_{\left\{ \left|(\mathbf{B}_1)_j\right| > \sqrt{K \log n}  \right\}} \prod_{\ell=1}^d \left(  1+ \left|(\mathbf{B}_1)_\ell \right|^{k_\ell} \right)\Big)\nonumber\\
	& \lesssim \Pi_\mathbf{0} \left( \left(  1+ \left|B_1\right|^J \right) e^{(1-\varepsilon)|B_1|^2/2} \right)
	\sum_{n=2}^\infty
	n^{ (2m- |\mathbf{k}|)/(4\kappa)} 	n^{-(1-\varepsilon)K/2}\lesssim
	\sum_{n=2}^\infty
	n^{ -(4\kappa +1)/(4\kappa)}<\infty.
\end{align}
Thus we have that $\P$-almost surely,
\begin{align}\label{step_22}
	&\lim_{n\to \infty}r_n^{m/2}
	\sup_{ r_n \leq s< r_{n+1}}\Big|e^{-(1+\frac{\Vert\theta\Vert^2}{2})\sqrt{r_n}} \sum_{u\in N(\sqrt{r_n}) } e^{-\theta \cdot \mathbf{X}_u(\sqrt{r_n})} 1_{\left\{ \cup_{j=1}^d\left\{ \left|(\mathbf{X}_u(\sqrt{r_n}))_j +\theta_j \sqrt{r_n}\right| > \sqrt{K\sqrt{r_n} \log n}  \right\}\right\}} \nonumber\\& \quad \quad \quad \times
	\prod_{j=1}^d \frac{\left(\sqrt{r_n}\right)^{k_j/2} }{s^{k_j/2}}   H_{k_j}\Big(\frac{(\mathbf{X}_u(\sqrt{r_n}))_j+\theta_j \sqrt{r_n}}{r_n^{1/4}}\Big)  \Big|=0.
\end{align}
Combining \eqref{step_21'} and  \eqref{step_22}, we get
\begin{align}
	&\mu_{s}^\theta \left((-\infty, \mathbf{b}\sqrt{s}]\right)  =  o(s^{-m/2}) +e^{-(1+\frac{\Vert\theta\Vert^2}{2})\sqrt{r_n}} \sum_{u\in N(\sqrt{r_n}) } e^{-\theta \cdot \mathbf{X}_u(\sqrt{r_n})} \nonumber\\
	& \quad \quad \times \prod_{j=1}^d\Big(  \sum_{k=0}^J \frac{(-1)^{k}}{k!} \frac{1}{s^{k/2}} \frac{\mathrm{d}^k}{\mathrm{d} b_j^k} \Phi(b_j) \Big(\left(\sqrt{r_n}\right)^{k/2} H_k\Big(\frac{(\mathbf{X}_u(\sqrt{r_n}))_j+\theta_j \sqrt{r_n}}{r_n^{1/4}}\Big)\Big)  \Big)
	\\
	&=  o(s^{-m/2}) + \sum_{\mathbf{k}: \mathbf{k}\leq J\mathbf{1}} \frac{(-1)^{|\mathbf{k}|}}{\mathbf{k}!} \frac{1}{s^{|\mathbf{k}|/2}} D^\mathbf{k} \Phi_d(\mathbf{b}) M_{\sqrt{r_n}}^{(\mathbf{k},\theta)}.
\end{align}
Since $\lambda> |J\mathbf{1}|/2 = d(3m+5)$, it follows from
Proposition \ref{Convergence-rate-Martingale} that
for any $\mathbf{k}\in \N^d$ with $m+1\leq |\mathbf{k}| \leq |J\mathbf{1}|$, $s^{-|\mathbf{k}|/2} M_{\sqrt{r_n}}^{(\mathbf{k},\theta)} = o(s^{-m/2})$. Thus
\begin{align}
	&\mu_{s}^\theta \left((-\infty, \mathbf{b}\sqrt{s}]\right)  =  \sum_{\mathbf{k}: |\mathbf{k}|\leq m} \frac{(-1)^{|\mathbf{k}|}}{\mathbf{k}!} \frac{1}{s^{|\mathbf{k}|/2}} D^\mathbf{k} \Phi_d(\mathbf{b}) M_{\sqrt{r_n}}^{(\mathbf{k},\theta)}+ o(s^{-m/2}).
\end{align}
Take $\eta>0$ sufficient small so that $\lambda>\frac{3m}{2}+\eta$.
Then by Proposition \ref{Convergence-rate-Martingale},
for any $\mathbf{k}\in \N^d$ with $0\leq |\mathbf{k}| \leq m$,
\begin{align}\label{Convergence-rate-2}
	M_{\sqrt{r_n}}^{(\mathbf{k},\theta)} - M_\infty^{(\mathbf{k},\theta)} = o(r_n^{-(\lambda -|\mathbf{k}|/2)/2+\eta/2}) = o\left(r_n^{-m/2} r_n ^{ m/2-(\lambda-m/2)/2 + \eta/2 }\right) = o(r_n^{-m/2}),
\end{align}
which implies that  as $s\to\infty$,
\begin{align}
	&\mu_{s}^\theta \left((-\infty, \mathbf{b}\sqrt{s}]\right)  =\sum_{\mathbf{k}: |\mathbf{k}|\leq m} \frac{(-1)^{|\mathbf{k}|}}{\mathbf{k}!} \frac{1}{s^{|\mathbf{k}|/2}} D^\mathbf{k} \Phi_d(\mathbf{b}) M_{\infty}^{(\mathbf{k},\theta)}+ o(s^{-m/2}),\quad
	\P\mbox{-a.s.}
\end{align}
Therefore,
the assertion of the theorem is valid under the assumption $\lambda> \max\left\{ 3m+8, d(3m+5)\right\}$.
\hfill$\Box$
\bigskip

\textbf{Proof of Theorem \ref{thm2}:}
Let $m\in \mathbb{N}$ and assume  \eqref{LlogL} holds for some $\lambda> \max\{d(3m+5), 3m+3d+8\}$.
Recall that $r_n=n^{1/\kappa}$ and  $\kappa=m+3$. Put $K:= m/\kappa +3$ and define
\[
Y_{s,n,u}:=   e^{-\theta \cdot  \mathbf{X}_u(\sqrt{r_n})}
\prod_{j=1}^d \frac{\sqrt{s}}{\sqrt{ s-\sqrt{r_n}}} \int_{a_j}^{b_j}  \phi\Big( \frac{z_j-\theta_j \sqrt{r_n}- (\mathbf{X}_u(\sqrt{r_n}))_j}{\sqrt{s-r_n}}\Big)\mathrm{d}z_j.
\]
Since $\lambda>3m+3d +8 $, by \eqref{step_24},
for
any $\theta\in \R^d$ with $\Vert\theta\Vert< \sqrt{2}$, any $\mathbf{a}, \mathbf{b}\in \mathbb{R}^d$ with
$\mathbf{a}<\mathbf{b}$ and
$s\in[r_n, r_{n+1})$,
$\mathbb{P}$-almost surely, as $s\to\infty$,
\begin{equation}\label{expan-mu1}
	s^{d/2}\mu_s^\theta \left( (\mathbf{a},\mathbf{b}]\right)=
	e^{-(1+\frac{\Vert\theta\Vert^2}{2})\sqrt{r_n}} \sum_{u\in N(\sqrt{r_n}) }  Y_{s,n,u}
	+ o(s^{-m/2}).
\end{equation}
Noticing that
$0\leq \sup_{ r_n \leq s< r_{n+1}} Y_{s,n,u}\leq e^{-\theta\cdot \mathbf{X}_u(\sqrt{r_n})} \prod_{j=1}^d \frac{(b_j-a_j)\sqrt{r_{n+1}}}{\sqrt{r_n -\sqrt{r_n}}} \lesssim e^{-\theta \cdot \mathbf{X}_u(\sqrt{r_n})}$,
and  using \eqref{step_26}, we get
\begin{align}
	&  \sum_{n=2}^\infty
	r_n^{m/2} \mathbb{E} \Big(	\sup_{ r_n \leq s< r_{n+1}} e^{-(1+\frac{\Vert\theta\Vert^2}{2})\sqrt{r_n}} \sum_{u\in N(\sqrt{r_n}) } Y_{s,n,u}1_{\left\{ \cup_{j=1}^d\left\{ \left|(\mathbf{X}_u(\sqrt{r_n}))_j +\theta_j \sqrt{r_n}\right| > \sqrt{K\sqrt{r_n} \log n}\right\} \right\}} \Big)\nonumber\\
	& \lesssim d
	\sum_{n=2}^\infty
	n^{m/(2\kappa)} \Pi_\mathbf{0} \left(|(\mathbf{B}_{\sqrt{r_n}})_1| > \sqrt{K \sqrt{r_n} \log n}\right) <\infty.
\end{align}
Therefore, $\mathbb{P}$-almost surely,
\begin{align}\label{step_27}
	s^{m/2}  e^{-(1+\frac{\Vert\theta\Vert^2}{2})\sqrt{r_n}} \sum_{u\in N(\sqrt{r_n}) } Y_{s,n,u}1_{\left\{ \cup_{j=1}^d\left\{ \left|(\mathbf{X}_u(\sqrt{r_n}))_j +\theta_j \sqrt{r_n}\right| > \sqrt{K\sqrt{r_n} \log n}\right\} \right\}} \stackrel{s\to\infty}{\longrightarrow} 0.
\end{align}
By \eqref{expan-mu1},  for $s\in [r_n, r_{n+1})$,
\begin{align}\label{expan-mu2}
	& s^{d/2}\mu_s^\theta \left( (\mathbf{a},\mathbf{b}]\right)\nonumber\\
	&= e^{-(1+\frac{\Vert\theta\Vert^2}{2})\sqrt{r_n}} \sum_{u\in N(\sqrt{r_n}) } Y_{s,n,u}1_{\left\{ \cap_{j=1}^d\left\{\left|(\mathbf{X}_u(\sqrt{r_n}))_j +\theta_j \sqrt{r_n}\right| \leq  \sqrt{K\sqrt{r_n} \log n}\right\} \right\}}+ o(s^{-m/2}).
\end{align}
Using Lemma \ref{lemma4},  on the set $\cap_{j=1}^d\{\left|(\mathbf{X}_u(\sqrt{r_n}))_j +\theta_j \sqrt{r_n}\right| \leq \sqrt{K\sqrt{r_n} \log n}\}$,
for $J=6m+10$,
\begin{align}
	&	\frac{\sqrt{s}}{\sqrt{ s-\sqrt{r_n}}} \int_{a_j}^{b_j}  \phi\Big( \frac{z_j-\theta_j \sqrt{r_n}- (\mathbf{X}_u(\sqrt{r_n}))_j}{\sqrt{s-r_n}}\Big)\mathrm{d}z_j \nonumber\\
	&= \int_{a_j}^{b_j} \phi\left(\frac{z_j}{\sqrt{s}}\right) \Big(\sum_{k=0}^J \frac{1}{k!} \frac{1}{s^{k/2}} H_k\Big(\frac{z_j}{\sqrt{s}}\Big)  r_n^{k/4}H_k\Big(\frac{ (\mathbf{X}_u(\sqrt{r_n}) )_j+\theta_j \sqrt{r_n}}{r_n^{1/4}}\Big) \Big) \mathrm{d} z_j+
	\varepsilon_{m,u,s,j}\nonumber\\
	& = \sum_{k=0}^J \frac{1}{k!} \frac{1}{s^{k/2}} \Big(\int_{a_j}^{b_j} \phi\Big(\frac{z_j}{\sqrt{s}}\Big)H_k\Big(\frac{z_j}{\sqrt{s}}\Big) \mathrm{d}z \Big) r_n^{k/4}H_k\Big(\frac{ (\mathbf{X}_u(\sqrt{r_n}))_j +\theta_j \sqrt{r_n}}{r_n^{1/4}}\Big) +
	\varepsilon_{m,u,s,j},
\end{align}
where
\[
r_{n+1}^{m/2}\sup_{s\in[r_n, r_{n+1})}\sup_{j\leq d}\sup_{u\in N(\sqrt{r_n})} |\varepsilon_{m,u,s,j} | 1_{\left\{ \cap_{j=1}^d\left\{\left|(\mathbf{X}_u(\sqrt{r_n}))_j +\theta_j \sqrt{r_n}\right| \leq  \sqrt{K\sqrt{r_n} \log n}\right\} \right\}}\to 0.
\]
Therefore, using \eqref{expan-mu2} and an argument similar to that leading to \eqref{step_21'}, we get
\begin{align}\label{expan-mu3}
	&s^{d/2}\mu_s^\theta (({\mathbf a}, {\mathbf b}])
	=  o(s^{-m/2}) + e^{-(1+\frac{\Vert\theta\Vert^2}{2})\sqrt{r_n}} \nonumber\\
	&\quad\quad\times \sum_{u\in N(\sqrt{r_n}) } e^{-\theta \cdot \mathbf{X}_u(\sqrt{r_n})}1 _{\left\{ \cap_{j=1}^d\left\{\left|(\mathbf{X}_u(\sqrt{r_n}))_j +\theta_j \sqrt{r_n}\right| \leq  \sqrt{K\sqrt{r_n} \log n}\right\} \right\}}\nonumber\\
	&\quad\quad \times \prod_{j=1}^d \Big\{  \sum_{k=0}^J \frac{1}{k!} \frac{1}{s^{k/2}} \Big(\int_{a_j}^{b_j} \phi\Big(\frac{z_j}{\sqrt{s}}\Big)H_k\Big(\frac{z_j}{\sqrt{s}}\Big) \mathrm{d}z_j \Big) r_n^{k/4}H_k\Big(\frac{ (\mathbf{X}_u(\sqrt{r_n}) )_j+\theta_j \sqrt{r_n}}{r_n^{1/4}}\Big)  \Big\}.
\end{align}
By Lemma \ref{General-many-to-one} and \eqref{step_28} and the fact that $ \left|\int_{a}^{b} \phi\left(\frac{z}{\sqrt{s}}\right)H_k\left(\frac{z}{\sqrt{s}}\right) \mathrm{d}z \right| \lesssim |b-a|$, we have that for any $\mathbf{k}\in \N^d$ with $0\leq |\mathbf{k}|\leq J$,
\begin{align}
	&	 \sum_{n=2}^\infty
	r_n^{m/2} \mathbb{E} \sup_{ r_n \leq s< r_{n+1}} \Big|e^{-(1+\frac{\Vert\theta\Vert^2}{2})\sqrt{r_n}} \sum_{u\in N(\sqrt{r_n}) } e^{-\theta \cdot \mathbf{X}_u(\sqrt{r_n})} 1_{\left\{ \cup_{j=1}^d \left\{ \left|(\mathbf{X}_u(\sqrt{r_n}) )_j+\theta_j \sqrt{r_n}\right| > \sqrt{K\sqrt{r_n} \log n}\right\}\right\}}  \nonumber\\& \quad \times
	\prod_{\ell=1}^d  \Big(\int_{a_\ell}^{b_\ell} \phi\Big(\frac{z_\ell}{\sqrt{s}}\Big)H_{k_\ell}\Big(\frac{z_\ell}{\sqrt{s}}\Big) \mathrm{d}z_\ell  \Big)  \frac{\left(\sqrt{r_n}\right)^{k_\ell/2} }{s^{k_\ell/2}}   H_{k_\ell}\Big(\frac{(\mathbf{X}_u(\sqrt{r_n}))_\ell+\theta_\ell \sqrt{r_n}}{r_n^{1/4}}\Big)  \Big|\nonumber\\
	&\lesssim \sum_{n=2}^\infty
	r_n^{m/2} r_n^{-|\mathbf{k}|/4}  \Pi_\mathbf{0} \Big( 1_{\{ \cup_{j=1}^d\{ |(\mathbf{B}_{\sqrt{r_n}})_j| > \sqrt{K\sqrt{r_n}\log n}\} \}} \prod_{\ell=1}^d \Big| H_{k_\ell}\Big(\frac{(\mathbf{B}_{\sqrt{r_n}})_\ell}{r_n^{1/4}}\Big)\Big| \Big) \nonumber\\
	& \lesssim \sum_{j=1}^d  \sum_{n=2}^\infty
	n^{(2m-|\mathbf{k}|)/(4\kappa)}
	\Pi_\mathbf{0} \Big( 1_{\{ |(\mathbf{B}_1)_j| > \sqrt{K \log n}\}} \prod_{\ell=1}^d\left(1+ |(\mathbf{B}_1)_\ell|^J\right)\Big)<\infty.
\end{align}
Thus $\mathbb{P}$-almost surely, as $s\to\infty$,
\begin{align}
	& \lim_{n\to \infty}r_n^{m/2}
	\sup_{ r_n \leq s< r_{n+1}} \Big|e^{-(1+\frac{\Vert\theta\Vert^2}{2})\sqrt{r_n}} \sum_{u\in N(\sqrt{r_n}) } e^{-\theta \cdot \mathbf{X}_u(\sqrt{r_n})} 1_{\left\{ \cup_{j=1}^d \left\{ \left|(\mathbf{X}_u(\sqrt{r_n}) )_j+\theta_j \sqrt{r_n}\right| > \sqrt{K\sqrt{r_n} \log n}\right\}\right\}}  \nonumber\\& \quad \times
	\prod_{\ell=1}^d  \Big(\int_{a_\ell}^{b_\ell} \phi\Big(\frac{z_\ell}{\sqrt{s}}\Big)H_{k_\ell}\Big(\frac{z_\ell}{\sqrt{s}}\Big) \mathrm{d}z_\ell  \Big)  \frac{\left(\sqrt{r_n}\right)^{k_\ell/2} }{s^{k_\ell/2}}   H_{k_\ell}\Big(\frac{(\mathbf{X}_u(\sqrt{r_n}))_\ell+\theta_\ell \sqrt{r_n}}{r_n^{1/4}}\Big)  \Big|=0.
\end{align}
Therefore, by \eqref{expan-mu3},  since $\lambda> 3m+3d+8$,
\begin{align}
	&s^{d/2}\mu_s^\theta ((\mathbf{a}, \mathbf{b}])
	=  o(s^{-m/2}) + e^{-(1+\frac{\Vert\theta\Vert^2}{2})\sqrt{r_n}} \sum_{u\in N(\sqrt{r_n}) } e^{-\theta \cdot \mathbf{X}_u(\sqrt{r_n})}\nonumber\\
	&\quad \times \prod_{j=1}^d \Big\{  \sum_{k=0}^J \frac{1}{k!} \frac{1}{s^{k/2}} \Big(\int_{a_j}^{b_j} \phi\Big(\frac{z_j}{\sqrt{s}}\Big)H_k\Big(\frac{z_j}{\sqrt{s}}\Big) \mathrm{d}z_j \Big) r_n^{k/4}H_k\Big(\frac{ (\mathbf{X}_u(\sqrt{r_n}) )_j+\theta_j \sqrt{r_n}}{r_n^{1/4}}\Big)  \Big\}\\
	& = o(s^{-m/2}) + \sum_{\mathbf{k}: \mathbf{k}\leq J\mathbf{1}} \prod_{j=1}^d \frac{1}{k_j!} \frac{1}{s^{k_j/2}} \Big(\int_{a_j}^{b_j} \phi\Big(\frac{z_j}{\sqrt{s}}\Big)H_{k_j}\Big(\frac{z_j}{\sqrt{s}}\Big) \mathrm{d}z_j \Big) M_{\sqrt{r_n}}^{(\mathbf{k},\theta)}.
\end{align}
Since $\lambda> \max\{d(3m+5), 3m+3d+8\}$,
using Proposition \ref{Convergence-rate-Martingale}
and argument similar to that used in the proof of Theorem \ref{thm1}, we get that
\begin{align}\label{expan-mu4}
	& s^{d/2}\mu_s^\theta ((\mathbf{a}, \mathbf{b}])  \nonumber\\
	&= o(s^{-m/2})+ \sum_{\mathbf{k}: |\mathbf{k}|\leq m} \prod_{j=1}^d \frac{1}{k_j!} \frac{1}{s^{k_j/2}} \Big(\int_{a_j}^{b_j} \phi\Big(\frac{z_j}{\sqrt{s}}\Big)H_{k_j}\Big(\frac{z_j}{\sqrt{s}}\Big) \mathrm{d}z_j \Big) M_{\infty}^{(\mathbf{k},\theta)}.
\end{align}
By Taylor's expansion, as $x\to 0$,
\begin{align}\label{Taylor-1}
	\phi(x)= \sum_{j=0}^m \frac{\phi^{(j)}(0)}{j!} x^j + o(x^m).
\end{align}
Note that $\phi^{(k)}(x)= (-1)^k H_k(x)\phi(x)$ and that, for each
$1\leq k\leq m$,
\begin{align}\label{Taylor-2}
	\phi(x)H_k(x) = (-1)^k \sum_{j=0}^{m} \frac{\phi^{(k+j)}(0)}{j!} x^j + o(x^{m}).
\end{align}
Combining \eqref{Taylor-1} and \eqref{Taylor-2},  we get
\begin{align}
	&\sum_{\mathbf{k}: |\mathbf{k}|\leq m} \prod_{j=1}^d \frac{1}{k_j!} \frac{1}{s^{k_j/2}} \Big(\int_{a_j}^{b_j} \phi\Big(\frac{z_j}{\sqrt{s}}\Big)H_{k_j}\Big(\frac{z_j}{\sqrt{s}}\Big) \mathrm{d}z_j \Big) M_{\infty}^{(\mathbf{k},\theta)}\nonumber\\
	& =  o(s^{-m/2}) + \sum_{\mathbf{k}: |\mathbf{k}|\leq m}  \frac{(-1)^{|\mathbf{k}|}M_{\infty}^{(\mathbf{k},\theta)}}{\mathbf{k}! s^{|\mathbf{k}|/2}}\int_{[\mathbf{a},\mathbf{b}]} \sum_{\mathbf{i}: |\mathbf{i}|\leq m} \frac{ \prod_{j=1}^d \phi^{(k_j+i_j)}(0)z_{j}^{i_j}}{\mathbf{i}! s^{|\mathbf{i}|/2}}\mathrm{d}z_1...\mathrm{d}z_d\\
	& = o(s^{-m/2}) + \sum_{\mathbf{k}: |\mathbf{k}|\leq m}  \frac{(-1)^{|\mathbf{k}|}M_{\infty}^{(\mathbf{k},\theta)}}{\mathbf{k}! s^{|\mathbf{k}|/2}}\int_{[\mathbf{a},\mathbf{b}]} \sum_{\mathbf{i}: |\mathbf{i} +\mathbf{k}|\leq m} \frac{ \prod_{j=1}^d \phi^{(k_j+i_j)}(0)z_{j}^{i_j}}{\mathbf{i}! s^{|\mathbf{i}|/2}}\mathrm{d}z_1...\mathrm{d}z_d.
\end{align}
Therefore, by \eqref{expan-mu4}, we conclude that
\begin{align}
	& s^{d/2}\mu_s^\theta ((\mathbf{a}, \mathbf{b}])  =
	o(s^{-m/2}) + \sum_{\mathbf{k}: |\mathbf{k}|\leq m}  \frac{(-1)^{|\mathbf{k}|}M_{\infty}^{(\mathbf{k},\theta)}}{\mathbf{k}! s^{|\mathbf{k}|/2}}\int_{[\mathbf{a},\mathbf{b}]} \sum_{\mathbf{i}: |\mathbf{i} +\mathbf{k}|\leq m} \frac{ \prod_{j=1}^d \phi^{(k_j+i_j)}(0)z_{j}^{i_j}}{\mathbf{i}! s^{|\mathbf{i}|/2}}\mathrm{d}z_1...\mathrm{d}z_d\nonumber\\
	& = \sum_{\ell=0}^m \frac{1}{s^{\ell/2}} \sum_{j=0}^\ell (-1)^j\sum_{\mathbf{k}: |\mathbf{k}|=j}\frac{M_{\infty}^{(\mathbf{k},\theta)}}{\mathbf{k}!} \sum_{\mathbf{i}: |\mathbf{i}|=\ell-j}\frac{ \prod_{j=1}^d \phi^{(k_j+i_j)}(0)}{\mathbf{i}!}\int_{[\mathbf{a},\mathbf{b}]} \prod_{j=1}^d z_{j}^{i_j}\mathrm{d}z_1...\mathrm{d}z_d.
\end{align}
\hfill$\Box$

\noindent

\begin{singlespace}

\end{singlespace}

\vskip 0.2truein
\vskip 0.2truein

\noindent{\bf Haojie Hou:}  School of Mathematical Sciences, Peking
University,   Beijing, 100871, P.R. China. Email: {\texttt
houhaojie@pku.edu.cn}

\smallskip

\noindent{\bf Yan-Xia Ren:} LMAM School of Mathematical Sciences \& Center for
Statistical Science, Peking
University,  Beijing, 100871, P.R. China. Email: {\texttt
yxren@math.pku.edu.cn}

\smallskip
\noindent {\bf Renming Song:} Department of Mathematics,
University of Illinois Urbana-Champaign,
Urbana, IL 61801, U.S.A.
Email: {\texttt rsong@illinois.edu}

\end{document}